\documentclass[12pt,a4paper]{article}
\usepackage{amsthm,amsmath,amsfonts}
\usepackage{latexsym}
\usepackage[T1]{fontenc}
\usepackage{a4wide}
\usepackage{times}
\usepackage{array}
\usepackage{hhline}
\numberwithin{equation}{section}

\input epsf

\theoremstyle{plain}
\newtheorem{theorem}{Theorem}[section]

\newtheorem{lemma}{Lemma}[section]
\newtheorem{corollary}{Corollary}[section]
\theoremstyle{definition}
\newtheorem{definition}{Definition}[section]

\def\d{\partial}
\def\al{\alpha}
\def\de{\delta}
\def\ep{\epsilon}
\def\be{\beta}
\def\vfi{\varphi}

\def\si{\sigma}
\newcommand{\bR}{\mathbb R}
\newcommand{\bN}{\mathbb N}

\def\cF{\mathcal F}
\def\cJ{\mathcal J}

\def\sgn{\hbox{\rm sgn}\,}

\def\zlo{\bar\circ\,}
\newcommand{\JR}{{\rm J}\mbox{-}{\rm R}}
\begin{document}
\title{Almost classical solutions to the total variation flow}
\author{Karolina Kielak, Piotr Bogus\l aw Mucha \& Piotr Rybka}

\maketitle

\begin{center}
 {\small Instytut Matematyki Stosowanej i Mechaniki, Uniwersytet Warszawski\\
ul. Banacha 2, 02-097 Warszawa, Poland\\
 E-mail: p.mucha@mimuw.edu.pl, p.rybka@mimuw.edu.pl}
\end{center}

\date{}

\bigskip

{\bf Abstract. } The paper examines one-dimensional total variation
flow equation with Dirichlet boundary conditions. Thanks to a new
concept of ``almost classical'' solutions we are able to determine
evolution of facets -- flat regions of solutions. A key element of our
approach is the natural regularity determined by nonlinear elliptic
operator, for which $x^2$ is an irregular function. Such a point of
view allows us to construct  solutions. We apply this idea to
implement our approach to numerical simulations for typical initial
data. Due to the nature  of Dirichlet data  any
monotone function is an equilibrium. We prove that each solution reaches
such steady state in a finite time. 

\medskip



\section{Introduction}
The equation which is the topic of this paper
\begin{equation}\label{i1}
u_t -\frac d{dx}\left(\sgn(u_x)\right) =0, \qquad u(a) =a_b,\quad u(b)=a_e.
\end{equation}
is a one-dimensional example of the total-variation flow.  The motivation to
study this problem is twofold: a) image analysis, see
\cite{rudin-etco}, \cite{med-book}, \cite{allard}; b) crystal growth 
problems, see \cite{angenent}, \cite{taylor}, \cite{kogigi},
\cite{margetis}. There are 
different physically relevant models,
where a similar to ours surface energy appears, but the
corresponding evolutionary problem is not necessarily set up, see
e.g. \cite{belik}.

Equation (\ref{i1}) may be  interpreted  as a steepest
descent of the total variation, i.e.  we can write
(\ref{i1}) as a gradient flow $u_t\in -\partial E(u)$ for a
functional $E$. This is why we can apply the abstract nonlinear
semigroup theory of Komura, see \cite{brezis},  \cite{barbu}, to
obtain existence of solutions. This has 
been performed by \cite{fukui-giga}, \cite{kogi}, \cite{kogigi} and
also by  \cite{med1}, \cite{med2}, \cite{med3}, \cite{med4}. However, the 
generality of this tool does not permit to study fine points of
solutions to (\ref{i1}). 

Solutions to (\ref{i1}) enjoy interesting properties, Fukui and Giga,
\cite{fukui-giga}, have noticed that facets persist. By a facet we
mean a flat part (i.e. affine) of the
solution with zero slope. Zero is exactly
the point of singularity of function $|\cdot|$.  This is why the
problem of facet evolution is not only nonlocal but highly
anisotropic. Our equation (\ref{i1}) is at least formally parabolic of
the second order. This is why we call the above  behavior of solutions 
the {\it sudden directional diffusion}. However, even more dramatic effects of singular diffussion can be seen in the fourth order problems, see \cite{giga-giga}

As we have already mentioned  some properties of facets were
established in \cite{fukui-giga}, e.g. their finite speed of propagation
was calculated. What is missing is the description of the process how
they merge and how they are created. In \cite{mury-non} we studied a
problem similar to (\ref{i1}). We worked there with a simplification of the
flow of a closed curve by the singular mean weighted curvature. 
%
We have shown existence of so-called almost classical
solutions, i.e. there is a finite number of time instances when the
time derivative does not exist.
However the results of \cite{mury-non} indicate lack of efficiency of
the methods used there. This fact is our motivation to rebuilt the
theory from the very beginning. For this reason we consider here the
model system admitting effects of sudden directional diffusion.  
Hoping that our approach will be suitable for more general systems.

Our approach is as follows. We notice that
the implicit time discretization leads to a series of Yosida approximations
to the operator on the right-hand-side (r.h.s. for short) of
(\ref{i1}). We  study them quite precisely, 
because we  consider variable time steps. As a result we 
capture the moment when two facets merge. We do not perform any
further special considerations.   
We want to see how the regularity of original solutions is transported via
solvability of the Yosida approximation. Due to the one-dimensional
character of the problem we are able to 
obtain  a result so good that it is of the maximal regularity
character, what is rather expected for quasilinear parabolic systems.  
Let us underline that properly understood smoothness 
is the most important question connected to solvability of the original
system. We have to modify standard  
regularity setting in order to capture all phenomena appearing in  the
system we study. As a result of our considerations we come to the
conclusion  that the best
smoothness we could expect for a solution $u$  that $u(\cdot, t)$ be
piecewise linear  
function, while $x^2$ is an example of an irregular function.

Our main goal is monitoring the evolution, as well creation, of the facets and a precise
description of the regularity of solutions to  (\ref{i1}), which we
construct here. For this purpose we apply methods, which are
distinctively different from those in the literature. We develop ideas
which appeared in our earlier works.  The key point is a construction
of a proper composition of two multivalued operators: the first one is
$\sgn$ understood as a maximal monotone graph, the other one is $u_x$,
which is defined only a.e. We leave aside the issue that in general
this is a measure, not a function. This problem is resolved differently
by the authors applying the semigroup approach,  \cite{fukui-giga},
\cite{med-book},  \cite{kogigi}, \cite{med4} etc. We treat  $u_x$ as a
Clarke differential (see (\ref{depaw}) and the text below this
formula).  Here,
we show that this composition is helpful when:\\
(a) we construct solutions, see Theorem \ref{lemYA};\\ 
(b) we discuss regularity of solutions, see Theorems \ref{th:main} and
\ref{th:main2}. \\
At the moment, however, the usefulness of this approach is limited to one
dimension. The advantage of our method is also simplicity, the
composition is explicitly computable. As an extra result we obtain asymptotics of solutions. The Dirichlet boundary conditions imply that the set of possible equilibria consists of monotone functions. Our analysis shows that steady state must be reached in finite time.

On the other hand, there are two sorts of results available up to now to deal with (\ref{i1}):\\
1) the method based on the abstract semigroup theory, see e.g.
\cite{fukui-giga}, \cite{med-book},  \cite{kogigi} and \cite{med4}.
It is  very general and elegant, it
enables us to 
study the facet motion, but it does not 
capture all relevant information. The intrinsic difficulty associated
with this method is the fact that the energy functional corresponding
to (\ref{i1}) is not coercive, also see below Lemma \ref{ixq} and the
proof of Theorem \ref{lemYA}.\\ 
2) the method based on the appropriate definition of the viscosity
solution \cite{giga-osher}. However, a different kind of problem was
studied there. This is an active research field, see \cite{giga-cp}.

Our approach is based on the Yosida approximation, defined by as a
solution of the resolvent problem 
\begin{equation}\label{i2}
\lambda u- \frac d{dx}\left(\sgn(u_x)\right)=\lambda v\qquad
\hbox{in }(a,b),\quad (u-v)|_{\partial[a,b]} =0.
\end{equation}

There are a couple of points to be made here. Firstly, we will
construct  $u$, a solution to  (\ref{i2}), by  very simple
means, this is done is Section \ref{secYo}. This process resembles 
looking for a good notion of a weak solutions to a PDE. Since we
came up with an integral equation we will call its solutions 
{\it mild} ones, see formula (\ref{na11}).

Secondly,  (\ref{i2}) may be interpreted as an Euler-Lagrange equation
for a non-standard variational functional. Namely, we set
\begin{equation}\label{i2a}
\cJ(u) =\left\{
\begin{array}{ll}
\int_a^b |Du| & \hbox{if }u\in D(\cJ)\equiv 
\{u\in BV[a,b],\ u(a) =a_b,\ u(b)=a_e\},\\
+\infty & \hbox{if } L_2(a,b)\setminus D(\cJ),
\end{array}\right.
\end{equation}
where $\int_a^b |Du| $ is the total variation of measure $Du$.  We
stress that we consider the space $BV$ over a closed interval.
Then,  (\ref{i2}) may be seen as
\begin{equation}\label{in1}
v \in u + h \partial \cJ(u),
\end{equation}
where $\partial \cJ$ is the subdifferential of $\cJ$ and
$h=\frac1\lambda$. We shall see 
that the well-established convex analysis will yield 
existence of a unique solution to inclusion (\ref{in1}). This solution will
be called {\it variational}. Since variational solutions are stronger
(we shall see this), thus both solutions coincide.

We note that
the Dirichlet problem in the multidimensional case is much more
difficult, in particular the meaning of the boundary condition is not
clear, see \cite{mazon2010}.

Thus, no matter which point of view we adopted, $u$ is given as
the action of the nonlinear resolvent operator $R(\lambda, A)$ on $v$, i.e.
$$
u = R(\lambda, A) (v) \equiv (\lambda  +   A)^{-1}(v),
$$  
where 
$A=-\frac\partial{\partial x} \sgn \frac\partial{\partial x}$. 
However, the notion of a mild solution to (\ref{i2}) does not permit
us to interpret this equation  easily. On the other hand, by convex
analysis, we can see (\ref{i2}) as an inclusion (\ref{in1}). 

The definition of the nonlinear resolvent operator leads to a detailed
study of  $\cJ$. 
One of our 
results, see Theorem \ref{wngl}
is a characterization of solutions to (\ref{i2}). The advantage of
(\ref{i2}) is that it permits to monitor closely 
behavior of  facets. It says that the regularity propagates. That
is, if $v$ is such that $v_x$ belongs to the $BV$  space and the number of
connected components of the properly understood set $\{x: v_x(x) =0 \}$ is finite, then $u_x$ has the same property for sufficiently
large $\lambda$.

It is well-known that the nonlinear resolvent leads to Yosida
approximation, which is the key object in the construction of the
nonlinear semigroup in the Komura theory. Namely, we set 
\begin{equation}\label{yosida}
A_\lambda u = \lambda (u - R(\lambda,A)(\lambda u) ).
\end{equation} 
Our observation is that a maximal monotone multivalued operator like
$\sgn$ taking 
values in $[-1,1]$ may be composed with a multifunction properly
generalizing a function of bounded total variation. We shall describe
here this composition denoted by $\bar\circ$, see Section \ref{komp}.
We introduced such an operation in \cite{mury-non}, see also
\cite{mury-cba}. We also point to an essential difficulty here, which is
the problem of composition of two multivalued operators. Even if
both of them are maximal monotone, the result need not be monotone
nor single valued. 
If the outer
of the two operators we compose is a subdifferential, then we expect
that the result is closely related to the minimal section of the
subdifferential. 

One of our main 
results says that $A_\lambda u$ defined by (\ref{yosida}) indeed converges
to $-\frac\partial{\partial x} 
\sgn \bar\circ\, u_x$. Moreover, we have an error estimate, see
Theorem \ref{lemYA}, formula (\ref{yerr}). In this way we justify
correctness of the 
new notion. Due to the  ``explicit'' nature of $\bar\circ$, we may
better describe  the regularity of solutions to (\ref{i2}). 

Once we have constructed the Yosida approximation, we show existence of
solution to the approximating problem $u^\lambda_t = -A_\lambda(u^\lambda)$ on
short time intervals, 
where $u^\lambda(t_0)$ is given. This is done in Lemma \ref{l:4.1}. 
In
fact, the method is close in spirit to the construction of the
nonlinear semigroup, see \cite{cran-lig}. 
Convergence of the approximate solutions is shown at the end of Section 4.
Here,
we use the full power of the Yosida approximation to capture the
finite number of time instances when the solution $u(t)$ is just right
differentiable with respect to time, otherwise the  derivative
exists. The point is
that we can control the distance to the original problem (\ref{i1}), so
that we can  monitor the time instances when facets merge.

Let us tell few words about the approach of proving our result. 
First, we define a space of admissible functions giving regularity of
constructed solutions. Furthermore, we state main results together with
an explanation of the meaning of almost classical solutions. In
Section 3 we study the Yosida approximation for our system,
concentrating on qualitative analysis of solutions. Proofs in this
part are based on a direct construction  which is possible due to the fine properties of chosen regularity. Subsequently, we prove the main results concerning existence and regularity. Section 5 is devoted to an alternative proof of existence for the resolvent operator based on the classical approach via Calculus of Variations.  This analysis shows that restrictions taken in Sections 3 and 4 are natural and reasonable.
Finally, we study the asymptotics of solutions and  present  an
example of an explicit solution. We conclude our paper with
 numerical simulations. They are based upon the
semidiscretization. Since they present a series of time snapshots,
these pictures contain only the round-off error. At each time step
there is no discretization error. The examples in  Section 6 present the typical behavior, for which each solution becomes a monotone function in finite time.

\section{The composition $\zlo$ and the main result }\label{komp}

Our main goal is to present a new approach to solvability of systems
of type (\ref{i1}). We construct a solution with this novel
technique, and next we compare it with ones obtained in a more
standard approach. This will 
clarify  why some assumptions, which seem to be artificial, after
deeper analysis will look completely natural. 
The total variation flow is a good example for such experiment, since
we know precisely the solution, additionally its simple form allows us to
deduce the  qualitative properties by standard methods of the
calculus of variation. 

The first step is to define the basic regularity class of functions.

\begin{definition}\label{def:a1} (cf. \cite[Chapter 5]{ziemer})
We say that a real valued function $u$, defined over a closed interval
$[a,b]$, {\it belongs to  $ BV[a,b]$}, provided that
$$
\|Du\| \equiv\int_a^b|Du| < \infty ,
$$
where $|Du|$ is the total variation of the measure $Du$. We recall
that
$$
\|u\|_{BV[a,b]} = \|Du\| + \| u\|_1.
$$
\end{definition}

\noindent
For the sake of definiteness, but without any loss of generality we
assume that $[a,b]=[0,1]$.

\medskip

Additionally, we  treat $BV$  functions as multi-valued
function. This is easy for functions which are derivatives, $u_x\in
BV[a,b]$. This is very useful in the regularity study of solution to
(\ref{i1}). Indeed, if $u$ and $u_x$ belong to $ BV[a,b]$, then $u$ is
Lipschitz continuous. Hence, $\frac{d^+u}{dx}$ and $\frac{d^-u}{dx}$
exist everywhere and they differ on at most countable set. Thus, we may set
\begin{equation}\label{depaw}
\partial_x u(s) = \{\tau u_x^- + (1-\tau) u_x^+:\ \tau\in [0,1]\}.
\end{equation}
Under our assumptions on $u$, the set $\partial_x u(x)$ is the Clarke
differential of $u$ and equality holds in (\ref{depaw}) due to
\cite[Section 2, Ex. 1]{chang}. If $u$ is convex, then $\partial_x u$
is the well-known subdifferential of $u$. As a result, if $u_x\in BV$,
then for each $ x_0\in (0,1)$, we have
$$
\partial_x u(x_0)=[\lim_{x\to x_0^-} u_x(x),\lim_{x\to x_0^+} u_x(x)]_{or},
$$
where $[a,b]_{or}=[a,b]$ for $a\leq b$ and $[a,b]_{or}=[b,a]$ for $b>a$.

However, the  description of solutions as functions whose derivatives
belong to $BV$  is not
sufficient. We have required to restrict our attention  to its subclass. 
There is a need to control the facets,
which we shall explain momentarily. {\it A facet} of $u$, $F$ is a
closed, connected piece
of  graph of $u$ with  zero slope,
i.e. $F=F(\xi^-,\xi^+)=\{(x,y):\ y=const=u(x_0),\ x\in[\xi^-,\xi^+]\}$, 
which is maximal with respect 
  to inclusion of sets. The interval $[\xi^-,\xi^+]$ will be called
  the {\it set of parameters} or {\it  preimage} of facet $F$.

Let us recall that zero is the only
point, where the absolute value, $|\cdot|$, the integrand in the
definition of $\cJ$, fails to be differentiable. Thus, the special role
of the zero slope and facets.

We shall also  distinguish a  subclass of facets. We shall say that a facet
$F(\xi^-,\xi^+)$ 
{\it has  zero curvature}, if and only if there is such $\ep>0$, that
function $u$ restricted 
to $[\xi^--\ep, \xi^++\ep]$ is monotone. In the case the function under
consideration is increasing this means that 
$u(\xi^--\epsilon)< u(\xi^-)=u(\xi^+)< u(\xi^++\epsilon)$. We shall
see that  zero curvature 
facets do not move at  all. There may be even an infinite number of
them. They have no influence on the evolution of  
the system. For that
reason we introduce the following objects, capturing the essential
phenomena. We shall say that a facet $F(\zeta^-,\zeta^+)$ of $u$
is an {\it essential facet}. It
will be denoted by $F_{ess}(\zeta^-,\zeta^+)$, provided that
there exists $\epsilon>0$ such that
either

\begin{quote}
$u$ is decreasing on $(\zeta^- - \epsilon,\zeta^-)$ and
$u(t)>u(\zeta^-)$ for $t \in (\zeta^- - \epsilon,\zeta^-)$ 
and $u$ is increasing on  $(\zeta^+,\zeta^+ + \epsilon)$ and
$u(t)>u(\zeta^+)$ for $t \in (\zeta^+,\zeta^+ + \epsilon)$ 
(then we call such a facet {\it convex}); moreover we set
\begin{equation}\label{tra-cx}
\sgn \kappa_{[\zeta^-, \zeta^+]} = 1
\end{equation}
\end{quote}

\noindent
or

\begin{quote}
$u$ is increasing on $(\zeta^- - \epsilon,\zeta^-)$ and
$u(t)<u(\zeta^-)$ for $t \in (\zeta^- - \epsilon,\zeta^-)$ 
and $u$ is decreasing on  $(\zeta^+,\zeta^+ + \epsilon)$ and
$u(t)<u(\zeta^+)$ for $t \in (\zeta^+,\zeta^+ + \epsilon)$ (then we
call such facet  {\it concave}); moreover, we set
\begin{equation}\label{tra-cc}
\sgn \kappa_{[\zeta^-, \zeta^+]} = -1.
\end{equation}
\end{quote}
It may happen that $\zeta^-=\zeta^+=:\zeta$, then we shall call
$F(\zeta,\zeta)$ a {\it  degenerate essential facet}. In this case
$u$ has a strict local minimum or a strict maximum at point $\zeta$.

We will call $\sgn \kappa_{[\zeta^-, \zeta^+]}$ the {\it transition number}
of facet $F(\zeta^-, \zeta^+)$. For the sake of consistency we set the
transition number $\sgn \kappa_{[\zeta^-, \zeta^+]}$ to zero for a zero
curvature facet $F(\zeta^-, \zeta^+)$. 
\smallskip

The union of parameter sets of all essential facets is denoted by  $\Xi_{ess}(w)$ and 
 $K_{ess}(w)$ is the number of essential facets, including degenerate facets.

\begin{definition}\label{def:a2} Let us suppose that $w=\partial_x u\in
BV[0,1]$, where $u$ is absolutely continuous and $\partial_x u$ is the
Clarke differential of $u$. We 
define $\Xi(w)=\{x\in [0,1] : 0  \in w(x)\}$.
We say that $w$ as above is {\it J-regular} or shorter $w \in \JR$ iff
the set $\Xi_{ess}(w)\subset \Xi(w)$ consists of a finite number of components,
i.e.
\begin{equation}\label{dwacztery}
 \Xi_{ess}(w)=[a_1,b_1] \cup \ldots \cup [a_{K_{ess}(w)},b_{K_{ess}(w)}] 
\quad\mbox{where } a_i \leq b_i
\end{equation}
and each interval  $[a^i,b^i]$ is an argument set of an essential (nondegenerate or degenerate) facet
$F(a^i,b^i)$. In particular, components of 
$\Xi(w)\setminus \Xi_{ess}(w)$ consists only of arguments of zero curvature
facets of $u$. 

Our definition in particular excludes functions with fast oscillations
like $x^2\sin \frac{1}{x}$.
We distinguished above a  subset of $BV$ functions.

Since degenerate facets will be treated as pathology, 
for given $w\in \JR$, we define
\begin{equation}\label{d1}
L(w)=\min\{ b-a :  [a,b] \mbox { is a connected component of }
\Xi_{ess}(w)  \}.
\end{equation}
Note that $L(w)=0$ iff there exists a degenerate facet of $u$.

The name J-regular refers to the regularity of the integrand in the
functional $\cJ$, which has singular point at $p=0$. $J$-regularity of
$w=\partial u_x$ means that function $u$ can be split into finite
number of subdomains where it is monotone. 

We also define the following quantity,
\begin{equation}
 \|w\|_{\JR[0,1]}=\|w\|_{BV[0,1]} +K_{ess}(w),
\end{equation}
where $K_{ess}(w)$ is the number of connected parts of $\Xi_{ess}(w)$, however,
this is not any  norm in this space.
\end{definition}


We start with the definition of a useful class of admissible functions.

\begin{definition}\label{dxx}  We shall say that a function $a$ is
{\it admissible}, for short  $a \in AF[0,1]$,   iff $a:[0,1]\to\bR$,
\begin{equation}\label{a1}
 \alpha=\partial_{x} a \mbox{ ~~with~~  } \alpha \in \JR \mbox{
   ~~and~~ } a(0)=a_b, a(1)=a_e.
\end{equation}
Here,  $\partial_x a$ denotes   the  set-valued Clarke
differential of $a$. 
\end{definition}

We note that the above definition restricts the behavior of admissible
function at the boundary of the domain. Namely, if $a\in AF$, then 
$a$ is monotone on an interval $[0,x_0)$ for some $x_0  \in (0,1)$ 
and either 
$$
a(x_0) > a(0) 
\mbox{~~~ or~~~ } 
a(x_0) < a(0) .
$$
By the same token, 
$a$ is monotone on an interval $(x_0,1]$ for some $x_0  \in (0,1)$
and either
$$
a(x_0) > a(1) 
\mbox{~~~ or~~~ }
a(x_0) < a(1).
$$


Thus,  the Dirichlet boundary condition makes immobile any facet
touching the boundary. Hence, such facets behave as if they had zero
curvature. 

\smallskip

A composition of multivalued operators requires proper preparations.
Due to the needs of our paper, we restrict ourselves to a definition of 
$$
\sgn \bar\circ\, \alpha
$$
for a suitable class of multivalued operators $\alpha$.
Of course, it is most important to define this composition in the
interior of the domain we work with. See also \cite{mury-non}, \cite{mury-cba}.

\begin{definition}\label{ddxx} 
Let us suppose that $a$ is admissible and $\partial_x \beta=\alpha \in
\JR[0,1]$. The definition of 
$\sgn \bar\circ\, \alpha$ is pointwise. Let us first consider $x\in
[0,1]\setminus \Xi_{ess}(\alpha)$.
Then, there exists an interval $(a,b)$ containing $x$ and such that
either $\beta$ is increasing on $(a,b)$ or decreasing. In the first case we set
\begin{equation}\label{a2}
   \sgn \bar\circ\, \alpha (x) = 1;
\end{equation}
if $\beta$ is decreasing on $(a,b)$, then we set
\begin{equation}\label{a3}
   \sgn \bar\circ\, \alpha (x) = - 1.
\end{equation}
We note that the set $[0,1]\setminus  \Xi_{ess}(\alpha)$ is a 
finite sum of open intervals, on each of them
function $\beta$ is monotone. Furthermore, the end points of $[0,1]$
can not belong to $ \Xi_{ess}(\alpha)$.

Now, let us consider $x\in \Xi_{ess}(\alpha)$, then there is 
$[p,q]$ a connected component of $\Xi_{ess}(\alpha)$ containing $x$. If 
$F(p,q)$ is a  convex facet of $\beta$, then we set,
\begin{equation}\label{a6}
 \sgn \bar\circ\, \alpha(x) = \frac{2}{q-p} x -\frac{2p}{q-p} -1
 \mbox{~~~~for~~} x \in [p,q].
\end{equation}
If $F(p,q)$ is a concave facet of $\al$, then  we set,
\begin{equation}\label{a7}
 \sgn \bar\circ\, \alpha(x) = -\frac{2}{q-p} x +\frac{2p}{q-p} +1
 \mbox{~~~~for~~} x \in [p,q] .
\end{equation}
\end{definition}

\medskip 

We have already mentioned that the Dirichlet boundary condition does not permit
any motion of the facet touching the boundary. Thus, effectively, they
behave like zero-curvature facets. Part 2 of Definition \ref{ddxx}
takes this into account. 

Now, we are in a position to state  main results being also a
justification of the notion of  almost classical solutions to our system. 

\begin{theorem}\label{th:main}
 Let $u_0\in AF[0,1]$, $L(u_{0,x})>0$ with $u_0(0)=a_b$ and $u_0(1)=a_e$,
 then the  system (\ref{i1}) admits unique solution in the sense
 specified by (\ref{na11}) and such that
\begin{equation}
 u_x \in L_\infty(0,T; \JR[0,1]).
\end{equation}
Moreover, $u$ is an \underline{almost classical solution}, 
i.e. it fulfills (\ref{i1}) in the following sense
\begin{equation}\label{a8}
 \begin{array}{lcr}
u_t-\d_x \sgn \zlo u_x=0 & in  & [0,1]\times (0,T),\\
u(0,t)=a_b, \qquad u(1,t)=a_e & for & t \in [0,T),\\
u|_{t=0}=u_0 & on & [0,1],
\end{array}
\end{equation}
where the time derivative in (\ref{a8}) exists for all time instances,
except for at most a finite number of exceptions, the $x$ derivative
exists for at most a finite number of exceptions.
Additionally, $u(\cdot,t)\in AF[0,1]$ for $t \in [0,T]$.
\end{theorem}

We study a second order parabolic equation with the goal of
establishing existence of almost classical solutions. This is why we
do not consider general data in $L_2$, but those which are more
natural for this problem, where the jumps in $u_x$ and their number
matter most. This is why we look for $u$, which not only belongs to
$BV$, i.e. $u(\cdot,t)\in BV$, but also $u(\cdot,t)\in AF$. In addition, the
necessity of introducing essential facets will be
explained.

An improvement of the above result, showing  a regularization effects, is the following

\begin{theorem}\label{th:main2}
 Let $u_0$ be as in previous Theorem above, but $L(u_{0,x})=0$. Then,
 there  exists a unique mild 
solution to (\ref{i1}), which is 
almost classical and it  fulfills (\ref{a8}). 
Furthermore, $L(u_x(t))>0$ for $t>0$. 
\end{theorem}

The second theorem shows that the class of functions with non-degenerate
facets is typical, and each initially degenerate essential facet
momentarily evolve into an nontrivial interval. Furthermore,
creation of such a singularity is impossible. In order to explain this phenomena let us analyze the following very important example related to analysis of nonlinear elliptic operator defined by
subdifferential of (\ref{i2a}). 

We first recall  the basic
definition. We say that $w\in \partial\cJ(u)$ 
iff $w \in L_2(a,b)$ and for all $h\in L_2(a,b)$ the inequality holds,
\begin{equation}\label{rndep}
\cJ(u+h) -\cJ(u) \ge (w,h)_2.
\end{equation}
Here $(f,g)_2$ stands for the regular inner product in $L_2(a,b)$. We also say
that $v\in D(\partial \cJ)$, i.e. $v$ belongs to the domain of $\partial \cJ$
iff $\partial \cJ(v)\neq\emptyset$.

We state here our  fundamental example. We recall (\ref{i2a}) and for
the sake of convenience we set $(a,b)=(-1,1)$. Then we make the
following observation. 

\begin{lemma}\label{ixq} Function
 $\frac 12 x^2$  does not belong to $D(\d \cJ)$.
\end{lemma}

{\it Proof.} If $\frac 12 x^2 \in D(\d \cJ)$, then there existed $w \in L_2(-1,1)$ such that for all $\phi \in C^\infty_0 (-1,1)$ and $t\in \bR$
\begin{equation}\label{do1}
 \int_{(-1,1)} (|x+t\phi_x|-|x|)dx \geq t\int_{(-1,1)} w \phi dx.
\end{equation}
We restrict  ourselves to $\phi$ such that 
\begin{equation*}
 \phi \in C^\infty_0(-\delta,\delta) 
\mbox{~~~~ and ~~~~ supp}\, \phi_x [-\delta,-\delta/2] \cup [\delta/2,\delta].
\end{equation*}
Additionally,
\begin{equation*}
 \phi_x(t) > 0\mbox{ ~ for ~ } t \in (-\delta,-\delta/2), 
\phi_x(t) <0 \mbox{ ~ for ~ } t \in (\delta/2,-\delta)
\end{equation*}
and
\begin{equation*}
 \phi(-\delta)=\phi(\delta)=0, \phi(t)=1 \mbox{ ~ for ~} t \in
 (-\delta/2,\delta/2). 
\end{equation*}
Next, let us observe that
\begin{equation}\label{do5}
 |x+t\phi_x(x)|-|x|=t \phi_x(x) \sgn x \mbox{ ~~~ for ~} |t\phi_x(x)|\leq \delta/2;
\end{equation}
we keep in mind that $\phi_x (t)=0 $ for $t \in (-\delta/2,\delta/2)$.

Thus, for such $\phi$ and $t$ the r.h.s. of (\ref{do1}) equals
\begin{equation}\label{do6}
 \int_{(-\delta/2,\delta/2)} (|x+t\phi_x(x)|-|x|) dx=\int_{(-\delta,-\delta/2)} t \phi_x \cdot (-1) dx+\int_{(\delta/2,\delta)}
t\phi_x\cdot (1) dx = -2t.
\end{equation}
Hence, we get
\begin{equation*}
 -2t \geq t \int_{(-\delta,\delta)} w \phi dx,
\end{equation*}
what implies for $t>0$ 
\begin{equation}\label{do8}
 2\leq -\int_{(-\delta,\delta)} w \phi dx \leq \int_{(-\delta,\delta)} |w| dx \to 0,
\end{equation}
since $w\in L_2(-1,1)$. Thus, we have just reached a contradiction. Hence,
$\frac 12 x^2$ can not belong to $D(\d \cJ)$. \qed

The full description of the domain of the subdifferential
$\partial\cJ$ of (\ref{i2a}) is beyond the scope of this paper. 
There is a description of $D(\cJ)$ for the multidimensional version of
the problem we consider, see e.g. \cite{mazon2010}. It is based on
Anzellotti's formula for integration by parts
\cite{anzellotti}. However, a direct characterization  of this set for
the one-dimensional problem seems to be missing even though
this functional has been studied in the 
literature.

At the end we mention a result describing  the asymptotics of solutions, proved in the last section.

\begin{theorem}\label{asymp} There is finite $t_{ext}>0$ such that the solution
  $u$ reaches a steady state at  $t_{ext}$, i.e. $u(t) = u(t_{ext})$
  for $t>t_{ext}$. Moreover, we have an explicit estimate for
  $t_{ext}$ in terms of $u_0$, see (\ref{exte}).
\end{theorem}

The above result shows that the limit of any solution, as
time goes to infinity, is always
a monotone function, and this will be proved and illustrated in Section 6. There
we present numerical simulations based on the analysis of system
(\ref{i1}). It is interesting to note that in comparison  with
\cite{feng} who deals with the multidimensional case, our computations
do not contain any discretization error. A rich possibility of stationary states are allowed thanks to Dirichlet boundary conditions.
Note that such picture is impossible for Neumann boundary constraints, for which there are only trivial/constant equilibria.

\section{Yosida approximation}\label{secYo}

The central object for our considerations is the Yosida approximation
to $-\partial_x\sgn \partial_x$. First, we introduce an auxiliary
notion of a nonlinear resolvent operator to the following problem,
\begin{equation}\label{y1}
 \lambda u-\frac{d}{dx}\sgn (u_{x})=\lambda v\mbox{ ~~~ on ~~ } [0,1],\qquad  u= v\quad\hbox{ at }\partial[0,1],
\end{equation}
where $v$ is a given element of $L_2(0,1)$.

\begin{definition}\label{d:1} An operator assigning to $v\in\JR$ a
unique  solution, $u\in\JR$,  to  (\ref{y1}) 
will be called {\it the resolvent of}
$A=-\partial_x\sgn \partial_x$ 
and we denote it  by    
$$
u=R(\lambda,A)v.
$$
\end{definition}


Now, we may  introduce the Yosida approximation to $A$.

\begin{definition}\label{d:2} Let us assume that 
$A=-\partial_x\sgn  \partial_x$  is as above and $\lambda>0$. 
An operator $A_\lambda:\JR \to \JR$ given by
$$
A_\lambda u=\lambda (u - R(\lambda,A)(\lambda u))
$$
is called the {\it Yosida approximation of $A$.}
\end{definition}

Since the notion of Yosida approximation seems well-understood, we will
use it to explain the meaning of $A$. For this purpose we will fix $w\in
\JR$ and $\lambda>0$. We set $u^\lambda:= R(\lambda,A)w$. We will look
more closely at $A_\lambda(u^\lambda)$.

\begin{theorem}\label{lemYA}
{\sl  Let us assume that $w\in AF[0,1]$, i.e. $w_x\in \JR$,  then there exists a unique solution to
\begin{equation}\label{d16}
\lambda u+A(u)=\lambda w\quad\hbox{in }(0,1),\qquad u(0)=w(0),\ u(1)=w(1),
\end{equation}
denoted by $u^\lambda$, fulfilling
\begin{equation}\label{y2}
 \|u^\lambda_{x}\|_{BV[0,1]} \leq   \|w_{x}\|_{BV[0,1]}.
\end{equation}
Moreover, there is $\lambda_0>0$ such that
$$
K_{ess}(u^\lambda_x)=K_{ess}(w_x) \mbox{~~~~~ for ~~} \lambda>\lambda_0 \mbox{ ~~  with ~~ } 
\|u^\lambda_x\|_{\JR}  \leq \|w_x\|_{\JR} .
$$
Furthermore, if $L(w_{x})=d>0$, equation (\ref{d16}) can be restated as follows
\begin{equation}\label{yerr}
\lambda u^\lambda-\partial_x \sgn \bar\circ u^\lambda_{x} =\lambda w +V(\lambda,x),
\end{equation}
where $V(\lambda,x) \to 0$ in $L_q$ for all $q<\infty$ as $\lambda
\to \infty$. In addition
$$
A_\lambda(u^\lambda)\to - \partial_x \sgn \bar \circ w_{x} \mbox{ ~~~ in ~~~ } L_q(0,1)
 \mbox{~~~with~~~} q<\infty.
$$
}
\end{theorem}
{\it Proof.} 
We would like to present an independent proof of existence of solutions to
system (\ref{d16}), which is based on simple tools, without any explicit
reference to calculus of variations. For this purpose, we restrict
ourselves to $w\in AF$ and for sufficiently large $\lambda$. 
A simple construction of $u^\lambda$ for a given $w$ based upon Lemma
\ref{podko}, is presented below.

Our assumptions give us 
\begin{equation}\label{na1}
 \Xi_{ess}(w_x)=\bigcup_{i=1}^{K_{ess}(w_x)} [a^i_*,b^i_*]
\end{equation}
with $a^i_* \leq b^i_*$.
Moreover, $a^1_*>0$ and $b_*^{K_{ess}(w_x)}<1$.

Below, we present a construction of $u^\lambda$. Namely, we
consider  system (\ref{d16}) in a neighborhood of preimage of an essential
facet $[a^i_*,b^i_*]$ of $w$ (it may be degenerate)  and we prescribe
the evolution of this facet. If $\lambda$ is  sufficiently large, then
we  keep the number $K_{ess}$ constant. 

\begin{lemma}\label{podko} Let us suppose that $w$ satisfies the
  assumptions of Theorem \ref{lemYA}. Then,  for sufficiently large
  $\lambda$, and for each $i=1,\ldots,K_{ess}(w_x)$  there exist
  monotone functions 
$$
\lambda\mapsto a_i(\lambda)\quad\hbox{ and }\quad
\lambda\mapsto b_i(\lambda),
$$ 
which are solutions to the following problem, 
\begin{equation}\label{na3}
 (b^i(\lambda)-a^i(\lambda))w(a^i(\lambda))=\int_{a^i(\lambda)}^{b^i(\lambda)}
w +2\frac {1}{\lambda} \sgn \kappa_{[a^i_*,b^i_*]}, \qquad
w(b^i(\lambda))=w(a^i(\lambda)). 
\end{equation}
These solutions are defined  locally, i.e. in a neighborhood of
$[a^i_*,b^i_*]$. 

We recall that, the transition numbers $\sgn \kappa_{[a^i_*,b^i_*]} $
were defined in 
(\ref{tra-cx}), (\ref{tra-cc}). Additionally, we require 
\begin{equation}\label{na4}
 a^1(\lambda) > 0, \qquad b^{K_{ess}(w_x)}(\lambda) < 1 \mbox{~~ and~~
 } b^i(\lambda) < a^{i+1}(\lambda)\quad\hbox{for }i=1,\ldots,K_{ess}(w_x)-1.
\end{equation}

However, if  $\lambda_0$ is the
greatest lower bound of $\lambda$ as above, then
one of the three possibilities occurs,
\begin{equation}\label{na5}
 a^1(\lambda_0)=0 \quad\mbox{or}\quad
 b^{K(w_x)}(\lambda_0)=1\quad\mbox{or}\quad
a^i(\lambda_0) = b^{i+1}(\lambda_0).
\end{equation}
It is worthwhile to underline that the lemma holds if $L(w_x)=0$, too.

\end{lemma}
{\it Proof.} 
Let fix $i$ in $\{1,\ldots,K_{ess}(w_x)\}$. Problem (\ref{na3}) comes from integration of equation (\ref{d16}) over a neighborhood of facet $[a^i_*,b^i_*]$. For $\tau\in \bR$ in a
neighborhood of zero and such
that $\tau\sgn\kappa_{[a_*^i,b_*^i]}>0$, we  set
\begin{equation}\label{potem}
\bar a^i(\tau) = \min (w|_{[b^{i-1}_*, a^i_*]})^{-1}(w(a_*^i)+\tau),\qquad
\bar b^i(\tau) = \max (w|_{[b^{i}_*, a^{i+1}_*]})^{-1}(w(b_*^i)+\tau).
\end{equation}
This definition is correct, because functions $w|_{[b^{i-1}_*, a^i_*]}$
and $w|_{[b^{i}_*, a^{i+1}_*]}$ are monotone. If these functions are
strictly monotone, then $w^{-1}(w(b_*^i)+\tau)$ is  strictly monotone
too, so the min/max are redundant. However, if there exists 
$\{\alpha\}\neq[\alpha,\beta]\subset \Xi (w)$ and $[\alpha,\beta]\subset
[b^{i-1}_*, a^i_*]$ (resp. $[\alpha,\beta]\subset [b^{i}_*,
a^{i+1}_*]$, then $(w|_{[b^{i-1}_*, a^i_*]})^{-1}$ (resp.
$(w|_{[b^{i}_*, a^{i+1}_*]})^{-1}$) is a maximal monotone graph and
min/max makes  
$\bar a^i(\cdot)$ (resp. $\bar b^i(\cdot)$) single valued and
discontinuous. However, the function
\begin{equation*}
\tau\mapsto (\bar b^i(\tau)-\bar a^i(\tau))w(\bar a^i(\tau))-
\int_{\bar a^i(\tau)}^{\bar b^i(\tau)}w (s)\, ds=:F_i(\tau), \qquad
i=1,\ldots,K_{ess}(w_x),
\end{equation*} 
is continuous. Indeed, if $\tau_0$ is  point, where $\bar a^i$ and $\bar b^i$
are continuous, then this statement is clear. Let us suppose that at
$\tau_0$ function $\bar a^i$ has a jump (the argument for $\bar b^i$
is the same). Then,  $[\bar a^i(\tau_0), \beta]\subset \Xi(w_x)$, where
$\bar a^i(\tau_0)<\beta$ and for any $x\in[\bar a^i(\tau_0), \beta]$ we
have
\begin{equation}\label{fsta}
(\bar b^i(\tau_0)-\bar a^i(\tau_0))w(\bar a^i(\tau_0))-
\int_{\bar a^i(\tau_0)}^{\bar b^i(\tau_0)}w (s)\, ds =
(\bar b^i(\tau_0)-x)w(x)-
\int_{x}^{\bar b^i(\tau)}w (s)\, ds.
\end{equation} 
This is so, because we notice that $w$ restricted to $[\bar
  a^i(\tau_0), \beta]$ is 
constant and equal to $w(a^i_*)+\tau_0$. Moreover,
$$
\int_{\bar a^i(\tau_0)}^{\bar b^i(\tau_0)}w (s)\, ds=
\int_{\bar a^i(\tau_0)}^x w (s)\, ds 
+\int_x^{\bar b^i(\tau_0)}w (s)\,ds =
(x-\bar a^i(\tau_0)) (w(a^i_*)+\tau_0)+\int_x^{\bar b^i(\tau_0)}w
  (s)\,ds.
$$
Hence, our claim follows, i.e. continuity of $F_i$, $i=1,\ldots,
K_{ess}(w)$. Indeed , let us suppose that $\tau_n$ converges from 
one side to $\tau_0$ (the side, left or right, depends upon
$\sgn \kappa_{[a^i_*,b^i_*]}$) so that $\lim_{n\to\infty} \bar a^i(\tau_n)
=\gamma$, where $\gamma=\bar a^i(\tau_0)$ or $\gamma= \beta$. Then,
due to (\ref{fsta}) we deduce continuity of $F_i$.


Subsequently, if we take $\lambda$ sufficiently large, then 
$\frac 2\lambda \sgn \kappa_{[a^i_*,b^i_*]}$ is in the range of $F_i$,
i.e. there exists $\tau_i=\tau_i(\lambda)$ such that
$F_i(\tau(\lambda)) =\frac 2\lambda \sgn \kappa_{[a^i_*,b^i_*]}$. If we
further make $\lambda$ larger, then we can make sure that
for each $i=1,\ldots,K_{ess}(w_x)$ we have 
$$
\bar b^{i-1}(\tau_i(\lambda))<\bar  a^i(\tau_i(\lambda))\quad\hbox{and}\quad
\bar b^{i}(\tau_i(\lambda))<\bar  a^{i+1}(\tau_i(\lambda)).
$$
Thus, we set
$$
a^i(\lambda) :=\bar  a^i(\tau_i(\lambda)),\qquad
b^{i}(\lambda):=\bar b^{i}(\tau_i(\lambda)).
$$
Let us define $\lambda_0$ to be the inf of $\lambda$'s constructed above.

We see that for $ \lambda_0$ one of the inequalities
$$
a^1(\lambda_0)>0,\qquad b^{i}(\lambda_0)< a^{i+1}(\lambda_0),\quad
i=1,\ldots,K_{ess}(w_x)-1,\qquad b^{K_{ess}(w_x)}(\lambda_0)< 1.
$$
become equality. $\Box$

\bigskip

This lemma permits us to define the function $u$ for $\lambda\ge\lambda_0$,
\begin{equation}\label{na6}
 u^\lambda=\left\{
\begin{array}{lcl}
\displaystyle w & \mbox{for}& x\in [0,1]\setminus \bigcup_{i=1}^{K_{ess}(w_x)}
[a^i(\lambda),b^i(\lambda)]\\
\displaystyle
w(a^i) & \mbox{for} & x \in [a^i(\lambda),b^i(\lambda)]
\end{array}
\right.
\end{equation}
We immediately notice that $K_{ess}(u^\lambda_x)=K_{ess}(w_x)$ and 
$\Xi_{ess}(u^\lambda_x)=\bigcup_{i=1}^{K_{ess}(u^\lambda_x)} [a^i,b^i]$, provided that 
$\lambda>\lambda_0$. 

Let us analyze what happens at $\lambda=\lambda_0$. We know that one
of the three possibilities in (\ref{na5}) occurs. 
We notice that if $a^1(\lambda_0)=0$ or $b^{K_{ess}(w_x)}(\lambda_0)=
1$, then  a facet  of $u^\lambda$  touches the
boundary. Subsequently this facet
becomes a zero curvature facet, for it is immobile. This is a simple
consequence of 
Dirichlet boundary conditions which do not admit any evolution of
facets touching the boundary. 

Let us look at the case
$b^i(\lambda_0)=a^{i+1}(\lambda_0)$ for an index $i$. Thus, we obtain
the phenomenon of facet merging. In both cases the structure of the
set $\Xi_{ess}(u^\lambda_x)$  will 
be different from $\Xi_{ess}(w_x)$. As a result,  we have
\begin{equation}\label{na7}
 K_{ess}(u^\lambda_x)< K_{ess}(w_x).
\end{equation}
It is worth stressing that at the moment $\lambda=\lambda_0$ more than
two facets may merge, so we can not control the decrease of number $K$. 
In this case we have to slightly modify (\ref{na6}), since the
structure of $\Xi_{ess}(u^\lambda_x)$ is different  from $\Xi_{ess}(w_x)$. It
is sufficient to notice that the number of elements in the
decomposition  (\ref{na1})  has decreased.
\smallskip

It is clear that  for $\lambda\ge \lambda_0$, we have
\begin{equation}\label{na8}
 K_{ess}(u^\lambda_x)\leq K_{ess}(w_x)
\end{equation}
and by the construction, (\ref{na6}) it is also obvious  that (see
Definition \ref{def:a1}) 
\begin{equation}\label{na9}
 \|D u^\lambda_x\|\leq \|D w_{x}\|.
\end{equation}
Note that the boundary conditions are given, so (\ref{na9}) controls
the whole norm of $u^\lambda$. 

Once we constructed a solution $u^\lambda$ by (\ref{na6}), we shall
discuss the question: in 
what sense does it satisfy  equation (\ref{i2}).  One hint is given in
the process of construction $a^i(\lambda)$ and  $b^i(\lambda)$. This is
closely related to ideas in \cite{MuRySIAM}. If we stick with
differential inclusions, then formula
\begin{equation}\label{na10}
 u-w-\frac{1}{\lambda} \frac{d}{dx} \sgn u_x \ni 0,
\end{equation}
leads to difficulties,
because we did not provide any  definition of the last term on the
left-hand-side (l.h.s. for short).

Here  comes our meaning of a {\it mild solution}: for each $x\in
[0,1]$, the following 
inclusion must hold
\begin{equation}\label{na11}
 \int_0^x (u-w)dx' - \left.  \frac {1}{\lambda} \sgn u_x\right|_0^x \ni 0.
\end{equation}
We shall keep in mind that at $x=0$, we have $u=w$ (for the sake of
simplicity of notation we shall suppress the superscript $\lambda$,
when this does not lead into confusion).

In order to show that $u$ fulfills (\ref{na11}), we will examine a
neighborhood of the first component  of $\Xi_{ess}(u_x)$, i.e.
$[a^1,b^1]$. We take $x\in[0,a^1)$, then $u=w$ on $[0,x]$. Thus, it is
  enough to  check whether $\frac 1\lambda (\sgn u_x(0) -\sgn u_x(x)) \ni
  0$. We notice 
  that on $[0,x]\subset[0,a_1)$ function $u$ is  monotone. As a result
    $\sgn u_x(0)$ and $\sgn u_x(x)$ may equal $1$ or $[-1,1]$, provided
    that $u$
    is increasing. If on the other hand, $u$ is  decreasing on
    $[0,x]$, then $\sgn u_x(0)$ and $\sgn u_x(x)$ are equal to $-1$ or
    $[-1,1]$. If any of these possibilities occurs, then (\ref{na11})
    is fulfilled. 

We shall continue, after assuming for the sake of definiteness  that
facet $F(a^1,b^1)$ is convex. The argument for a concave facet is analogous.

Let us consider $x\in[a^1,b^1]$. We  interpret $\sgn t$ as
a multivalued function such that $\sgn 0=[-1,1]$. Then, we have for
$x\in [a^1,b^1]$ 
\begin{equation}\label{na13} 
 \int_0^x (u-w)dx'-\frac{1}{\lambda}[-1,1]+\frac{1}{\lambda} \sgn u_x|_{x'=0} \ni 0.
\end{equation}
Since we assumed that
the facet $F(a^1,b^1)$ is
convex, from (\ref{na3}) we find that
\begin{equation}\label{na14}
 0\leq \int_0^x (u-w)dx' \leq \frac 2\lambda.
\end{equation}
By the assumption we know that $\sgn u_x|_{x'=0} \ni -1$. Hence,
\begin{equation}\label{na15}
 \int_0^x (u-w)dx' - \frac {1}{\lambda} \in \frac {1}{\lambda} [-1,1].
\end{equation}
This shows (\ref{na11}) again. In case $F(a^1,b^1)$ is
concave, the argument is analogous.

 Let us now consider  $x\in (b^1,a^2]$,  then we have
\begin{eqnarray}\label{na12}
\int_0^x (u-w)dx' - \left.  \frac {1}{\lambda} \sgn u_x\right|_0^x &=& 
\int_0^{a^1} (u-w)dx' - \left.  \frac {1}{\lambda} \sgn u_x\right|_0^{a^1} +
 \int_{a^1}^{b^1} (u-w)dx' \nonumber\\
&&- \left.  \frac {1}{\lambda} \sgn u_x\right|_{a^1}^{b^1} +
 \int_{b^1}^x (u-w)dx' 
- \left.  \frac {1}{\lambda} \sgn u_x\right|_{b^1}^x\\&= 
 &I_1+ I_2+ I_3 .\nonumber
\end{eqnarray}
Here, we do have the freedom of choosing $\sgn u_x$ at $x=b^1$. Namely
we set $\sgn u_x(b^1) =-1$. We also know that $\sgn u_x(a^1) =1$ .

We recall that by the very construction of $a^1$ and $b^1$, we have
$I_2=0$. Subsequently, we notice that the argument performed for
$x\in[0,a^1)$ applies also to $x\in(b^1,a^2]$, Thus,
\begin{eqnarray*}
I_1+ I_2+ I_3  &=& - \frac1\lambda(1-\sgn u_x(0)) + 0
 - \frac1\lambda(-1+\sgn u_x(x))\\
 &=& \frac1\lambda(-\sgn u_x(0)+\sgn u_x(x)) \ni 0,
\end{eqnarray*}
i.e. (\ref{na11}) holds again. 

Repeating  the above procedure for each subsequent facet, we prove that
$u$ given by (\ref{na8}) fulfills (\ref{na11}).  
The case $x\in[b^{K_{ess}},1]$ is handled in the same way. Thus, we
proved the first part of Theorem \ref{lemYA} concerning existence.

We shall look more closely at the solutions when  $\lambda=\lambda_0$. We have
then  two basic possibilities:\\
(1) The first facet $F(a^1,b^1)$ or the last one $F(a^k,b^k)$  touches
the boundary, i.e. $a^1=0$ or resp. $b^k=1$. If this happens, then
$F(0,b^1)$, resp. $F(a^k,1)$, has zero curvature.\\
(2) Two or more facets merge, i.e. there are $i,r>0$ such that
$$
\lim_{\lambda\to\lambda_0} b^{i-1}(\lambda)= b^{i-1}(\lambda_0)<
a^i(\lambda_0)= \lim_{\lambda\to\lambda_0} a^i(\lambda)
$$ 
and
$$
\lim_{\lambda\to\lambda_0} b^{i+j}(\lambda)=
\lim_{\lambda\to\lambda_0} a^{i+1+j}(\lambda),\quad j=0,1,\ldots, r-1,
$$
and
$$
\lim_{\lambda\to\lambda_0} b^{i+r-1}(\lambda)<
\lim_{\lambda\to\lambda_0} a^{i+r}(\lambda).
$$
We adopt the convention that $b^0=0$ and $a^{k+1}=1$.

When this happens, we have two further sub-options:\\
(i) an odd number of facets merge, then
$F(a^i(\lambda_0),b^{i+r}(\lambda_0))$ has zero curvature;\\
(ii) an even number of facets merge, then
$[a^i(\lambda_0),b^{i+r}(\lambda_0)]\subset  \Xi_{ess}(u_x)$.

\noindent
Of course, it may happen that simultaneously a number of events of
type (2i) or (2ii) occurs.

First let us observe  that $u=w$  away from the set
$\{u_x=0\}$, so  we conclude  
$\Xi(w_x) \subset  \Xi(u_x)$. More precisely, the
equality holds on a larger set. Namely, if $F(a^i,b^i)$ is a zero
curvature facet and $\lambda>\lambda_0$, then the very construction of
$a^i(\lambda)$, $b^i(\lambda)$ implies that $u=w$ on $[a^i,b^i]$. 
If $u(a)=u(b)=w(a)=w(b)$, so there must be a point $c\in (a,b)$ such that
$0\in w_x$. Thus, we obtain  for any $\lambda>0$ 
$$
K_{ess}(u_x)\leq K_{ess}(w_x).
$$

\bigskip

Let $L(w_x)=d>0$, then we consider 
\begin{equation}\label{d0}
\lambda u + A_\lambda(u)=\lambda w \mbox{ ~~~~for $ \lambda>\lambda_0$},
\end{equation}
where we suppressed the superscript $\lambda$ over $u$.

As we have already seen taking large $\lambda$, i.e. $\lambda>\lambda_0$,  
excludes the possibility of facet merging or hitting  the boundary, thus 
$K_{ess}(w_x)=K_{ess}(u_x)$. Let us emphasize that $K_{ess}(u_x)$ may 
decrease only a finite number of times.

Let us suppose that $[a^*,b^*]$ is a connected component  of
$\Xi_{ess}(u_x)$, i.e. $a^*= a^{i_0}(\lambda)$, $b^*=
b^{i_0}(\lambda)$ for an index $i_0$. Without loss 
of generality, we may assume that this facet is 
convex. So, integrating (\ref{d0}) over $[a^*,b^*]$, we find
\begin{equation}\label{d2}
\int_{a^*}^{b^*} \lambda u - \int_{a^*}^{b^*} \lambda w =2.
\end{equation}

First, we want to find an answer to the following question. What we can
say about the behavior of the following quantity 
$\int_{a^*}^a+\int_b^{b^*} (\lambda u-\lambda w)$,
where $[a,b]$ is a connected component  of $\Xi_{ess}(w_x)$ contained
in $[a^*,b^*]$. In fact we assume, that $a= a^{i_0}$, $b=b^{i_0}$.

Since $d=L(w_x)$ is fixed and positive we find from (\ref{d2}) that
$$
2=\int_{a^*}^{b^*} \lambda (u-w) \geq \int_a^b \lambda (u-w) \geq d
\lambda (u-w)|_{[a,b]}, 
$$
Because $u-w$ is monotone on $[a,b]$. 
As a result,
\begin{equation}\label{d4}
\lambda (u-w)|_{[a,b]}\leq \frac 2d.
\end{equation}
%
Then we conclude that
\begin{equation*}
\int_b^{b^*} \lambda (u-w) \leq (b^*-b)\lambda[w(b^*)-w(b)].
\end{equation*}
At the same time (\ref{d4}) yields, $w(b^*)-w(b) \leq
\frac{2}{d\lambda }$. On the other hand, $w$ is  monotone on set
$(b,b^*)$. Hence (\ref{d4}) implies that
\begin{equation}\label{d6}
b^*-b \equiv b^{i_0}(\lambda)-b^{i_0} \leq W^{-1}(\frac{2}{d\lambda}),
\end{equation}
where $W^{-1}(\cdot)$ is a strictly monotone (possibly multivalued) function, equal 
$w^{-1}$ (restricted to an interval of monotonicity) plus a constant such that $\lim_{t\to 0^+}W^{-1}(t)=0$.
Eventually, we get
\begin{equation}\label{d7}
\int_b^{b^*} \lambda (u-w) \leq W^{-1}(\frac{2}{d\lambda}) \frac{2}{d} \to 0 \mbox{~~as~~} \lambda\to \infty.
\end{equation}
Since the analysis for $(a^*,a)$ is the same, hence (\ref{d6}) and
(\ref{d7}) imply that
\begin{equation*}
\int_{a^*}^{b^*} \lambda (u-w) = 2 +V(\lambda ) 
\mbox{ ~~~ with ~~~ } V(\lambda ) \to 0 \mbox{~~as~~} \lambda\to \infty.
\end{equation*}
Note that  $V(\lambda )$ depends only on $w$, so in Section
\ref{secflo} we will study the approximation  error $V(\lambda )$ and
we will show uniform bounds, provided that $L(w_x) \geq d >0$. 

Integrating (\ref{d0}) yields
\begin{equation}\label{d9}
\int_{a^*}^{b^*} \lambda (u-w)=\int_{a^*}^{b^*} - A_\lambda(u) = 2,
\end{equation}
but the pointwise information from the equation yields
\begin{equation}\label{d10}
\lambda(u-w)|_{[a,b]}=-A_\lambda(u)= const.
\end{equation}
Thus, taking into account (\ref{d9}) and (\ref{d10}), we get
\begin{equation*}
- A_\lambda(u)|_{[a,b]} \to 2/(b-a) \qquad a.e. \mbox{~~as~~} \lambda\to \infty.
\end{equation*}
Here, we used that $a^*= a^{i_0}(\lambda)\to a^{i_0}$, $b^*=
b^{i_0}(\lambda)\to b^{i_0}$ as $\lambda$ goes to infinity.  But $a,b$
depends only on $w$, additionally we shall keep 
 in mind that (\ref{d4}) via (\ref{d0}) implies that $|A_\lambda(u)| \leq 2/d$ on whole $[0,1]$.

Clearly, by Definition \ref{ddxx}
$$
\partial_x \sgn \bar\circ u_{x} = \frac{2}{(b^*-a^*)} \mbox{ ~~ for ~~ } x \in [a^*,b^*].
$$
Hence, we have proved that
\begin{equation}\label{d11}
A_\lambda(u)=- \partial_x \sgn \bar\circ u_{x} +V(\lambda,x),
\end{equation}
where $V(\lambda) = \int_{a^*}^{b^*}V(\lambda,x)\,dx$ and
$V(\lambda,x) \to 0$ in at least $L_1(I)$. Here, we should note
clearly 
that all depend on $\lambda$, since $a^*=a^{i_0}(\lambda)$, 
$b^*=b^{i_0}(\lambda)$. We see that we have already proved
that 
$|V(\lambda,x)| \leq \frac 2d,$ and $ \mu(\{\hbox{supp}\, V(\lambda,.)\}) \to 0$
which gives a relatively strong convergence. 
Note that in (\ref{d11}) we are not able to obtain  ``pure''
discontinuity in the composition $\bar\circ$,  since we work with
solutions only,  hence $\sgn \bar \circ u_x^\lambda$ must be piecewise linear.

Next question is: whether
$\partial_x \sgn \bar\circ u_x^\lambda \to - \partial_x \sgn \bar\circ w_x$
and in which space?

Let us observe that (see Definition \ref{def:a1})
\begin{equation}\label{d13}
 \|Du^\lambda_x\| \leq \|D w_x\| \mbox{ ~~~ and ~~~ } u^\lambda_{x} \to w_{x} \mbox{ in measure on } I.
\end{equation}
It follows that
\begin{equation*}
  \sgn \bar \circ u^\lambda_{x} |_{\Xi(w_x)} \to  \sgn \bar \circ
  w_{x}|_{\Xi(w_x)} \mbox{ uniformly}. 
\end{equation*}
We remember that $\sgn \bar \circ u^\lambda_{x}$ and $\sgn \bar \circ
w_{x}$ are piecewise linear functions and the set $\Xi(w_{x})$ is
independent from $\lambda$,
but the case $L(w_x)=d>0$ implies that 
\begin{equation}\label{d15}
 A(u^\lambda) \to - \frac{d}{dx} \sgn \bar\circ w_x \mbox{ in }
 L_q(0,1)\qquad q\in[1,\infty).
\end{equation}
Theorem \ref{lemYA} is proved.
\qed

In particular, as a result of our analysis, we get that the constructed
solution to (\ref{d16}) is variational. 

\begin{lemma}\label{l:var}
Function $u^\lambda$ given by Theorem \ref{lemYA} is a variational solution to (\ref{d16}), i.e. $u^\lambda$ fulfills
\begin{equation}\label{d15a}
 (\lambda u^\lambda,\phi)+(\sigma(x),\phi')=(\lambda w,\phi) \mbox{ ~~ for each ~ } \phi \in C^\infty_0(0,1)
\end{equation}
and $\sigma(x) \in \sgn \circ u_x(x)$, where here $\circ$ denotes the standard composition.
\end{lemma}

{\it Proof.} From the inclusion (\ref{na11}), we are able to find such
$\sigma$ that 
\begin{equation}
 \int_0^x (u-w) -\frac{1}{\lambda} \sigma (x) + \frac{1}{\lambda} \sigma (0)=0.
\end{equation}
Then, testing it by $\phi'$ with $\phi \in C_0^\infty(0,1)$, we get
(\ref{d15a}). In particular, we already have shown that $\lambda
R(\lambda,A)\lambda$ is a monotone operator in $L_2$. \qed 

\section{The construction of the flow}\label{secflo}
A key point of our construction of solution is an approximation of the
original problem based on the Yosida approximation. Here, we meet techniques characteristic for the homogeneous Boltzmann equation \cite{diBlasio,mucha}.  For given $\lambda$,
$t_0$ and $A_\lambda$ defined in (\ref{yosida}), we introduce the following
equation for $u^\lambda$,
\begin{equation}\label{c1}
 u^\lambda(t+t_0)=u^\lambda(t_0)-\int_{t_0}^{t_0+t} A_\lambda(u^\lambda) \,ds,
\quad u^\lambda(0,t_0+t)=a_b, \quad u^\lambda(1,t_0+t)=a_e \mbox{ for } t \in (0,T).
\end{equation}
We stress that its solvability, established below, does note require
that  $L(u_x(t_0))>0$.  

\begin{lemma}\label{l:4.1}
Let us suppose that $u^\lambda(t_0)\in \JR(I)$, where $I=[0,1]$, then
there exists a unique 
solution $u^\lambda$ to (\ref{c1}) on the time interval 
$(t_0,t_0+\frac{1}{3\lambda})$ and
$$
u^\lambda \in C(t_0,t_0+ \frac{1}{3\lambda}; L_2(I)) .
$$
Moreover,
\begin{equation}\label{c2}
 \sup_{t\in (0,\frac{1}{3\lambda})} \|u^\lambda(t_0+t)\|_{\JR} \leq \|u^\lambda(t_0)\|_{\JR}.
\end{equation}
\end{lemma}

{\it Proof.} We will first show the bounds. Let us suppose that $u^\lambda$ is
a solution to (\ref{c1}), then 
Definition \ref{d:2} and  the observation
$\frac{d}{dt}[e^{\lambda t}u^\lambda]= - e^{\lambda t}A_\lambda(u^\lambda)+\lambda e^{\lambda t}u^\lambda$ imply that,
\begin{equation}\label{c3}
 u^\lambda(t_0+t)=e^{-\lambda t} u^\lambda(t_0) + 
\int_{t_0}^{t_0+t} e^{-\lambda (t_0+t-s)} \lambda R(\lambda,A)\lambda u^\lambda(s) ds. 
\end{equation}
In order to obtain the estimate in $BV$, we apply Theorem \ref{lemYA},
inequality (\ref{y2}), getting
\begin{eqnarray}
 \sup_t  \|u^\lambda_{x}\|_{BV} 
&\leq &e^{-\lambda t} \|u^\lambda_{x}(t_0)\|_{BV}
+ \sup_t \|R(\lambda,A)\lambda u^\lambda(t)\|_{BV}
\int_0^t \lambda e^{-\lambda s}ds  \nonumber
\\[10pt] &
\leq &
e^{-\lambda t} \|u^\lambda_{x}(t_0)\|_{BV} 
+\sup_t \frac {1}{\lambda} \|\lambda
u^\lambda_x(t)\|_{BV}(1-e^{-\lambda t}).
\nonumber
\end{eqnarray}
So we get
\begin{equation}\label{c4a}
\sup_t \|u^\lambda_{x}\|_{BV} \leq \|u^\lambda_{x}(t_0)\|_{BV}.
\end{equation}

In order to prove existence, we fix $\lambda$ (we will omit the index $\lambda$ in
the considerations below) and we define a map 
$\Phi: C(0,T; L_2(I))\to C(0,T;L_2(I))$ 
such that
 $v=\Phi(w)$, 
where
\begin{equation}\label{c6}
 v(t)=e^{-\lambda t} v_0 + \int_0^t e^{\lambda (t-s)} \lambda R(\lambda,A)\lambda w ds.
\end{equation}
We notice that due to $\Xi((\lambda R(\lambda,A)\lambda w)_x)\supset
\Xi(w_x)$ we obtain 
$\Xi(v_{0,x}) \subset \Xi(w_x(t))$ for $t\in (0,T)$, provided that $w|_{t=t_0}=v_0$. Combining this observation
with $w|_{t=t_0}=v_0$ again yields,
\begin{equation}\label{c6a}
 \Xi(v_{0,x})\subset \Xi(v_x(t)) \mbox{ ~~ for ~~  } t \in (0,T).
\end{equation}

We see that a fixed point of the above map yields a solution to
(\ref{c1}) after a shift of time. For the purpose of proving existence
of a fixed point of $\Phi$, we will check that
$\Phi$ is a contraction. We notice that if $w$, $\bar w\in
C(0,T; L_2(I))$, then monotonicity of $R(\lambda,A)\lambda$
(thanks to Lemma \ref{l:var}) implies that 
$$
\| R(\lambda,A)\lambda w  - R(\lambda,A)\lambda\bar w\|_{L_2} \le \|w-\bar w\|_{L_2}.
$$
Hence, 
\begin{equation*}
\begin{array}{ll}
 \| \Phi(w) - \Phi(\bar w)\|_{L_\infty(0,T;L_2(I))} 
&\leq \int_0^t \lambda e^{-\lambda (t-s)}ds \|R(\lambda,A)\lambda w -R(\lambda,A)\lambda\bar
 w\|_{L_\infty(0,T;L_2(I))} 
\\[8pt]
&\leq ( 1- e^{-T\lambda}) \|w-\bar w\|_{L_\infty(0,T;L_2(I))},
\end{array}
\end{equation*}
i.e. $\Phi$ is a contraction
provided that $0<T \leq \frac{1}{3\lambda}$. Now, Banach
fixed point theorem  implies 
immediately  existence of $u^\lambda$, a unique solution to (\ref{c1})  in
$C(0,T;L_2(I))$.

An aspect is that
the solution to (\ref{c3}) can be recovered as a limit of the
following iterative process
\begin{equation}\label{c8}
v^{k+1}=\Phi(v^k).
\end{equation}

We have to show that the fixed point belongs to a better space. For
this purpose we use estimate (\ref{c4a}), which shows also that if
$\|v^0_{x}\|_{BV}=M$, then $\|v^k_{x}\|_{BV} \leq M$ for all $k \in
\bN$. Moreover, convergence in $L_2(I)$ implies convergence in
$L_1(I)$ and lower semicontinuity of the total variation measure (see
\cite[Theorem 5.2.1.]{ziemer}) yields $u^\lambda\in L_\infty(0,T; BV(I))$.

Finally we  show that 
\begin{equation}\label{c9}
K_{ess}(u^\lambda(t_0+\frac{1}{3\lambda}))\leq K_{ess}(u(t_0)).
\end{equation}
For this purpose  it is enough to prove that
$$
u^\lambda(t_0+t)=u^\lambda(t_0) \mbox{ ~~ on ~~} 
I \setminus  \Xi(u^\lambda(t_0+t)) \mbox{ for all } t\leq \frac{1}{3\lambda},
$$
but Theorem \ref{lemYA} implies 
$$
R(\lambda,A)\lambda u^\lambda =\lambda u^\lambda \mbox{~~~ on ~~} 
I\setminus \Xi(R(\lambda,A)\lambda u^\lambda), 
$$
namely $A_\lambda(u^\lambda)=0$ at 
$I\setminus \Xi(R(\lambda,A)\lambda u^\lambda)$. Additionally
(\ref{c6a}) yields that  
$\Xi(u^\lambda(t_0)) \subset \Xi(u^\lambda(t_0+\frac{1}{3\lambda})$, 
what finishes the proof of (\ref{c9}).

Thus, the definition of the solution to (\ref{c1}) as the limit of the
sequence (\ref{c8}) together with (\ref{c9}) imply (\ref{c2}). The
Lemma  is proved. \qed

\begin{lemma}\label{l:4.2} Let us consider $u^\lambda(\cdot)$ given by Lemma
  \ref{l:4.1}. If $L(u^\lambda(t_0))=0$, then  $L(u^\lambda(t_0+\frac{1}{3\lambda}))>0$. 
\end{lemma}

{\it Proof.} Let us assume a contrary, then
there exists a degenerate facet $F[a^i,b^i]$ with $a^i=b^i$ such that
all functions $u^\lambda(t_0+t)$ are convex in a neighborhood 
$(p,q)$ of point $a^i$ and they all have a minimum only in point $a^i$. If
functions $u^\lambda(t_0+t)$ are  concave, then the argument is
analogous. Let us then integrate (\ref{c1}) over $(a',b')$ such that
$a^i \in (a',b') \subset (p,q)$,
\begin{equation*}
 \int_{a'}^{b'} u^\lambda(t_0+t) =\int_{a'}^{b'} u^\lambda(t_0)
 -\int_{t_0}^{t_0+t} \int_{a'}^{b'} A_\lambda(u^\lambda) ds. 
\end{equation*}
But 
\begin{equation*}
 \int_{a'}^{b'} A_\lambda(u^\lambda)=\int_{a'}^{b'} \lambda (u^\lambda - R(\lambda,A)\lambda u^\lambda)=-2,
\end{equation*}
because $u^\lambda$  is convex on $(a',b')$. Hence, we find
\begin{equation*}
  \int_{a'}^{b'} u^\lambda(t_0+t) =\int_{a'}^{b'} u^\lambda(t_0) + 2t.
\end{equation*}
But if our assumption that $a^i=b^i$ were true, then we would be
allowed  to pass to the limits, $a'\to {a^i}^-$ and $b' \to {a^i}^-$
concluding that   
$0=0+2t$, which is impossible for positive $t$. Thus, we
showed that $u^\lambda(t_0+\frac{1}{3\lambda})$ does not admit
degenerate facets. \qed 

\bigskip
After these preparations, we finish
{\it the proofs of Theorems \ref{th:main} and \ref{th:main2}.}
We shall construct an approximation of solution on a fixed  time
interval, say $[0,1]$.
Let us assume that
$$
U^\lambda:[0,1]\times I \to \bR
$$
is given as follows
$$
U^\lambda=u^\lambda_k \mbox{ ~~ for ~~ } t \in [\frac{k}{3\lambda},\frac{k+1}{3\lambda}) \mbox{ ~~ and ~~ } 0\leq k < 3\lambda,
$$
where functions $\{u^\lambda_k\}$ are given by the following relations
$$
u^\lambda_1(t)=u_0 - \int_0^t A_\lambda(u^\lambda_1)ds \mbox{~~~~for~~} t \in (0,\frac{1}{3\lambda}],\\
$$
$$
u^\lambda_2(t_1+t)=u_1(t_1) - \int_{t_1}^{t_1+t} A_\lambda(u^\lambda_2)ds \mbox{~~~~for~~} t \in (0,\frac{1}{3\lambda}],\\
$$
$$
...
$$
$$
u^\lambda_{k+1}(t_k+t)=u_k(t_k) - \int_{t_k}^{t_k+t} A_\lambda(u^\lambda_{k+1})ds \mbox{~~~~for~~} t \in (0,\frac{1}{3\lambda}],\\
$$
$$
...
$$
$$
u^\lambda_{3\lambda}(t_{3\lambda-1}+t)=u^\lambda_{3\lambda-1}(t_{3\lambda-1}) - \int_{t_{3\lambda-1}}^{t_{3\lambda-1}+t}  A_\lambda(u^\lambda_{3\lambda})ds \mbox{~~~~for~~} t \in (0,\frac{1}{3\lambda}]
$$
and  $t_k=\frac{k}{3\lambda}$ for $0\leq k <3\lambda$.
$$
\| U^\lambda\|_{L_\infty(0,T;\JR)} \leq \|u_0\|_{\JR}.
$$

Now, we  pass to the limit with $\lambda$. The estimates  imply that $\|U^\lambda\|_{L_\infty(0,T;L_2(I))}\le M$. Thus, we
can extract a subsequence such that
$$
U^\lambda \rightharpoonup^* U^* \mbox{ weakly $*$ in } L_\infty(0,1;L_2(I)) .
$$
Moreover, the lower semicontinuity of the total variation measure
yields
$$
\|U^\lambda(t)\|_{BV} \le \|u(0)\|_{BV}\qquad\hbox{~~for a.e.~} t\in[0,1].
$$
Thus, we should look closer at the limit  
$$
U^*(t_0+t)=U^*(t_0) -  \lim_{\lambda\to \infty} \int_{t_0}^{t_0+t} A_\lambda(U^\lambda(t_0+t))\, ds.
$$


Let us observe that for a fixed $\lambda$ the function $K_{ess}(U^\lambda(t))$,
taking values in $\mathbb N$, is decreasing, so facet merging  may
occur just  only a finite number of times. 


Let $K(u_0)=k^0$, then for a given $\lambda$ we define $T^\lambda_1$ as follows
\begin{equation}\label{cc1}
K_{ess}(U^\lambda(t))=k^0 \mbox{ ~~~ for ~~~ } t \in [0,T^\lambda_1) \mbox{ ~ and ~ } K_{ess}(U^\lambda(T^\lambda_1))< k^0.
\end{equation}
For a subsequence
$\lim  T^\lambda_1 =:T_1$. Indeed $T^\lambda_1=T^{\lambda'}_1$ for all sufficiently
large $\lambda,\lambda'$ see Lemma \ref{lac4}, so we have here $T_1 >0$. However, we prefer to consider a
more general argument valid for more complex operators.

In a similar manner to (\ref{cc1}) we define a sequence of time instances
$\{T_k\}_{k=1}^{m}$. By the definitions, for any $\epsilon>0$ there
exists $\lambda_\epsilon$, such that for $\lambda>\lambda_\epsilon$ -- up to possible subsequence -- 
we can split the time interval $[0,1]$ into following parts
$$
[0,1)=[0,T_1-\epsilon)\cup[T_1-\epsilon,T_2+\epsilon)\cup [T_2+\epsilon,T_3-\epsilon)]\cup ... \cup [T_m+\epsilon,1)
$$
and 
$$
K_{ess}(u^\lambda(t))=K_{ess}(U^*(t)) \mbox{ for } t \in [T_k+\epsilon,T_{k+1}-\epsilon),
$$
so $\{T_k\}$ is a finite sequence of moments of time at which  facets
merge. In order to avoid unnecessary problems we restrict ourselves to
a suitable subsequence guaranteeing the above properties. 

Now, Theorem \ref{lemYA}  yields
$A_\lambda(U^\lambda) \to A(U^*)= - \d_x\sgn \zlo U^*_{x}  \mbox{ in } L_q(0,1)$ on time intervals $(T_k+\epsilon,T_{k+1}-\epsilon)$, 
since by (\ref{d11}) we control this convergence uniformly at whole intervals.
So we get
$$
U^*(t_0+t)=U^*(t_0) - \int_{t_0}^{t_0+t} A(U^*(s))ds,
$$
because we consider one interval $[T_k+\epsilon,T_{k+1}-\epsilon)$. However,
  crossing $T_k$ requires some extra care.

In order to extend the result on the whole interval $[0,1]$, it is
sufficient to prolong the solution onto interval
$[T_k-\epsilon,T_k+\epsilon)$. For this purpose we can use  that
  $u^\lambda$ belongs to 
  $C(0,1;L_1(I))$, see Lemma \ref{l:4.1}. Continuity of of the
  solution  allows us to cross points
  $T_k$. 
It follows that 
$$
\frac{d}{dt} U^* \mbox{ exists except points } \{T_k\}
$$
and by the properties of solutions on intervals $[T_k,T_{k+1})$ 
we find that the right-hand-side time derivative exists everywhere, including points $\{T_k\}$
$$
\frac{d}{dt^+} U^* \mbox{ exists everywhere on } [0,1].
$$
Finally, we have shown that $U^*$ fulfills
\begin{equation}
 \frac{d}{dt^+} U^*= - \frac{d}{dx} \sgn \zlo U^*
\end{equation}
as an almost classical solution. 

By construction $U^*(t) \in AF$, additionally Lemma \ref{l:4.2} yields $L(U^*(t))>0$ for $t>0$, even as $L(u_{0,x}) =0$. 
Moreover, the features of almost classical  solutions imply that they are variational, too. Hence, the monotonicity of $\sgn$ 
implies immediately uniqueness to our problem.
Theorems \ref{th:main} and \ref{th:main2} are proved. \qed

Now we want to obtain the same result starting from the classical point of view of the calculus of variation in order to explain the chosen regularity.

\section{The variational problem}\label{secvar}

In this section we will prove Theorem \ref{lemYA} using the
tools of the Calculus of Variations. This result establishes
existence of solutions to (\ref{y1}), i.e.
\begin{equation*}
 \lambda u-\frac{d}{dx}\sgn (u_{x})=\lambda v\hbox{ in } (0,1),\qquad 
 u=v\hbox{ for } x=0,1 
\end{equation*}
for an appropriate $v$.

Some parts of the argument, when $v\in\JR$ with $L(v_x)>0$ are
a repetition of results from Section \ref{secYo}. 
However, this repetition is necessary in order to explain that
approach from previous sections are
based on a reasonable class of function, which can be viewed as typical.

It is clear that first we have to give meaning to this equation. We
can easily see that it is formally an Euler-Lagrange equation for a
functional $\cJ_{h,v}$ defined below.
$$
\cJ_{h,v}(u) = h\cJ(u)+ \frac12\int_a^b (u-v)^2,
$$
where  $\cJ$ is introduced in (\ref{i2a}). 
When no ambiguity arises, we shall write $\cJ_v$ in place of
$\cJ_{h,v}$.

We notice that $\cJ_v$ is proper and convex. Momentarily, we shall see
that it is also lower semicontinuous, 
hence its subdifferential is well defined, see \cite{brezis}, in
particular $D(\partial\cJ)\neq\emptyset$. We recall that $u\in
D(\partial\cJ)$ if and only if $\partial\cJ(u)\neq\emptyset$.
It is a well known fact that $u$ is a minimizer of $\cJ_v$ iff 
\begin{equation}\label{rn1}
h\partial\cJ(u) + u-v \ni 0.
\end{equation}
Since $h\partial\cJ(\cdot) + Id$ is maximal monotone, then for any
$v\in L_2$
there exists $u\in D(\partial\cJ)$ satisfying (\ref{rn1}), see
\cite{brezis}.

In this way, we obtain our first interpretation of   (\ref{y1}) as a
differential inclusion. This is not very satisfactory as long as we do
not have a description of the regularity of the elements of
$D(\partial\cJ)$. We note the
basic observation and present its direct proof.

\begin{lemma}\label{bla}{\sl (a) For any $v\in L_2(a,b)$ functional
    $\cJ_v$ is lower semicontinuous in $L_2$.\\
(b) If $v\in L_2(a,b)$, then there exists
    $u\in D(\cJ)\subset BV(a,b)$ a unique minimizer of $\cJ_v$. Moreover,}
$$
\|Du\|=\int_0^1 |D u| \le  |B-A| + \frac 1{2h} \int_0^1 (v-\ell)^2 dx,
$$
where $\ell$ is an affine function such that $\ell(a) =A$, $\ell(b)
=B$.
\end{lemma}

{\it Proof.} (a) Let us suppose that $\{u_n\}\subset L_2 $ is a
sequence converging to $u$ in $L_2$. If 
$$
\liminf_{n\to\infty} \|Du_n\|
=\infty,
$$ then there is nothing to prove. Let us suppose then that 
$
\sup_{n\in \bN}\|Du_n\| \le K.
$
By the lower semicontinuity of the $BV$ seminorm, we infer that $u\in BV$
and $ \|Du\| \le K$. The problem is to
show that the limit $u$ satisfies the boundary conditions. 

If $v\in BV[a,b]$, then there is a
representative such that $\|D\tilde v\| = V_a^b(\tilde v)$. Moreover,
$\hbox{ess\,sup}\,| v|$ is finite, see \cite[Chapter 5]{ziemer}. Thus,
there is a 
representative $\bar v$ satisfying the boundary conditions and
$V_a^b(\bar v)\le \|D v\| + 4 \| v\|_{\infty} $. As a result, we 
select a sequence of representatives $\bar u_n$ satisfying the boundary
conditions and with uniformly bounded variations. Since $\bar u_n$
is  a sequence of bounded functions with commonly bounded total 
variation we  use Helly's theorem to deduce existence of subsequence
$\{u_{n_k}\}$ which converges to $u^\infty$ everywhere. Since all
functions $\{u_{n_k}\}$  satisfy the boundary data, the pointwise
limit will satisfy them too. Moreover, due to uniqueness of the limit
$u^\infty = u$ a.e. thus we can select a representative belonging to
$D(\cJ)$ as desired.  

(b)
By definition $\cJ_v$ is bounded below. Let us suppose
that $\{u_n\}$ is a minimizing sequence in $L^2$. Of course $u_n$'s
belong to $BV(a,b)$ and
$$
\int_0^1 |D u_n| + \frac 12 \int_0^1(u_n-v)^2 dx \le  K.
$$
i.e. the sequence  $\{u_n\}$ is bounded in the $BV$ norm. Since sets 
bounded in  $BV$ are compact in any $L_p(0,1)$, $p<\infty$, see
\cite{Attouch},  we deduce
existence of a subsequence $\{u_{n_k}\}$ converging to
$u$. Because of part (a)
we infer that $u\in D(\cJ)$ and
$$
\int_0^1 |D u| \le \liminf_{k\to\infty} \int_0^1 |D u_{n_k}|
\le K.
$$
Combining this with strong convergence of $\{u_{n_k}\}$ in $L_2$ we
come to the conclusion that  $u$ is a minimizer of $\cJ_v$. 

Uniqueness of a minimizer is a result of strict monotonicity of the
operator $Id +h\partial\cJ$. 

Since, $u$ is a minimizer, then $\cJ_v(\ell)\ge \cJ_v(u)$, where
$\ell$ is an affine function such that $\ell(a) =A$, $\ell(b)
=B$. Hence, the desired estimate follows.
\qed

We shall establish how much of the smoothness of $v$ is passed to
$u$. 
Here is our first observation.

\begin{theorem}\label{twac}
{\sl If $v\in W^{1}_{p}(a,b)$, where $p\in(1,\infty)$, then $u$
the unique minimizer of 
$$
\cJ_v(u)\equiv
\int_a^b h|u_x| + \frac12 (u-v)^2 \equiv h\cJ(u)+\int_a^b \frac12 (u-v)^2
$$
belongs to $W^{1}_{p}$ and $\|u\|_{1,p}\le \|v\|_{1,p}$.}
\end{theorem}

We want to look at the propagation of regularity, so the assumption
$v_x\in BV$ is natural from many possible view points. So here is our
main result, it will be shown after Theorem \ref{twac}. Its proof
follows from the analysis of the argument leading to  Theorem \ref{twac}.
 
\begin{theorem}\label{wngl}  {\sl Let us suppose that $v \in AC[a,b]$ and
$u$ be the corresponding minimizer of $\cJ_v$. Then,\\
(a) $K_{ess}(u)\le K_{ess}(v)$; \\
(b) if $v_x\in BV$ and $ K_{ess}(v)$ is finite, then  $u_x\in BV$ and
$
\|u_x\|_{BV} \le \|v_x\|_{BV}.
$
}
\end{theorem}

\bigskip
We see from its statement that a type of regularity which propagates
is defined by $v_x\in BV$ and a finiteness of the number  $K_{ess}(v)$. At
this point, we do not claim that this is optimal.

In order to provide a proof of Theorem \ref{twac}, we will proceed in
several steps. First we shall deal 
with continuous piecewise smooth functions, then we shall show that
our claim is true for any $v$ which may be approximated in $W^{1}_{2}$
by such functions. We need a simple  device to check that a function
is indeed a minimizer.

\begin{lemma}\label{lac0}{\sl Let us suppose that $v,u\in AC[a,b]$
    with $v(a)=u(a)$, $v(b)=u(b)$ and
    there exists $\si \in W^1_1(a,b)$ and such that $\si(x)\in
    \sgn(u_x(x))$ with $\sgn$ understood as a multivalued graph, which
    satisfies the equation 
\begin{equation}\label{rn}
h\frac d{dx} \si = u-v
\end{equation}
in the $L_1$ sense. Then, $u$ is a minimizer of $\cJ_v$.}
\end{lemma}
{\it Proof.} Let us take any $\vfi\in C^\infty_0$. Let us calculate
\begin{eqnarray*}
\cJ_v(u+\vfi) - \cJ_v(u) & = & 
h\int_a^b|u_x+\vfi_x| - h\int_a^b|u_x| +
\int_a^b\frac12 [(u+\vfi-v)^2 -(u-v)^2] \\
&\ge& h\int_a^b|u_x+\vfi_x| - h\int_a^b|u_x| +
\int_a^b (u-v)\vfi\\
 & = & h\int_a^b|u_x+\vfi_x| - h\int_a^b|u_x| -h\int_a^b \si \frac
d{dx}\vfi.
\end{eqnarray*}
We used (\ref{rn}) and the integration by parts. We deal separately
with the sets $\{u_x>0\}$, $\{u_x<0\}$ and $\{u_x=0\}$. We have,
\begin{eqnarray*}
(\cJ_v(u+\vfi) - \cJ_v(u))h^{-1} &\ge&
\int_{\{u_x>0\}}(|u_x+\vfi_x| - u_x -1\cdot \vfi_x) \\&+&
\int_{\{u_x<0\}}(|u_x+\vfi_x| + u_x +1\cdot \vfi_x) +
\int_{\{u_x=0\}}(|\vfi_x| - \si\cdot \vfi_x) \ge 0.
\end{eqnarray*}
We used here the fact that $\si(x)\in[-1,1]$ as well. 

Now, we deal with general $\vfi\in BV$ such that $u+\vfi\in D(\cJ)$. We
proceed by smooth approximation $\vfi_n$ such that $\vfi_n$ converges
to $\vfi$ in $L_1$ and $\| D\vfi_n\| \to \| D\vfi\|$. By what we have
already shown, we have
$$
\cJ_v(u+\vfi_n) \ge \cJ_v(u).
$$
Hence, the inequality is preserved after a passage to the limit.
Our claim follows. \qed

\bigskip
We may now start the regularity analysis.

\begin{lemma}\label{lac1}{\sl Let us suppose that $v\in C[a,b]$,
    $v(a)=A$, $v(b)=B$, and its derivative exists almost everywhere
    and it
    is piecewise continuous, its one sided derivatives exist
    everywhere and the sets $\{v_x>0\}$, $\{v_x<0\}$ are open and 
the number of essential facets of $v$ is finite.
Then, for any positive $h$  and $u$ a unique minimizer of
    $\cJ_{h,v}$, we have $u\in W^{1}_{p}$ with
\begin{equation*}
\| u\|_{1,p}\le \|v\|_{1,p}.
\end{equation*}
Moreover, there exists $\si\in W^{1}_{\infty}$, such that for all $x\in[a,b]$ we
have $\si(x)\in \sgn(u_x(x))$ and equation
(\ref{rn})
is satisfied everywhere except a finite number of points. In
addition,} 
$$
\|\si \|_{1,\infty} \le 1 + \frac 1h \|v\|_\infty.
$$
\end{lemma}
{\it Proof. } We shall proceed by induction. We first show, however, a
slightly stronger result if $v$ is monotone i.e. the number $K_{ess}$
is zero, and to fix attention we
assume that it is increasing. Namely, we
claim that in this case $u=v$. We have to show that for any $\vfi$ such
that $v+\vfi\in  D(\cJ)$, i.e. $\vfi$ is zero at the ends of $[0,1]$, we have
$$
\cJ_{h,v}(v+\vfi) \ge \cJ_{h,v}(v).
$$
Let us notice that  
\begin{eqnarray*}
\cJ_{h,v}(v+\vfi)&=&\int_a^b (h|v_x+\vfi_x| + \frac12 \vfi^2)
\ge \int_a^b h(v_x+\vfi_x)\\ 
&=& B-A  = \int_a^bh v_x =\cJ_{h,v}(v).
\end{eqnarray*}
We may also set $\si =1,$ since $v$ is increasing.

The first non trivial case occurs when we have a single essential
facet $F_{ess}(a,b)$. The 
set $[0,1]\setminus [a,b]$ consists of exactly two components
$E^+(v)$ and $E^-(v)$. They are such that $v|_{E^+(v)}$ is
increasing while $v|_{E^-(v)}$ is decreasing.
We stress that the endpoints $0$, $1$ cannot belong to any essential
facets. 
For the sake of fixing attention, we may assume that for all
$x_0\in[a,b]$ function $v$ has a maximum at $x_0$, $v_M=
\max v(x) = v(x_0)$. 
We can find $\xi^-\in E^-(v)$, $\xi^+\in E^+(v)$, i.e. $v$
increasing on $[\xi^-,a]$ while it is decreasing on $[b,\xi^+$], and
such that
\begin{equation}\label{vcom}
v(\xi^-) = v(\xi^+)=v_{com}
\end{equation} and 
$v_{com}$ is the smallest number with this property. In addition,
since $v$ is not strictly monotone on $E^+(v)$ or $E^-(v)$, we
require that if $\zeta\in E^+(v)$ (respectively, $\zeta\in E^-(v)$)
is another number satisfying (\ref{vcom}), then $\zeta\le \xi^+$
(respectively, $\zeta\ge \xi^-$). In this way $\xi^+$, $\xi^-$ are
uniquely defined.

We want to solve (\ref{rn}), for this purpose we will utilize results
of Lemma \ref{podko}.  In the present case the term $-\frac
{d\sigma}{dx}$ is used in place of $A(u)$. Since we are dealing with a
single facet we may be more specific about the range of $\tau$
appearing  (\ref{potem}).
We notice that for any $\tau\in (0,v_M-v_{com}]$ there exist
$\xi^-(\tau)\in[\xi^-,a]$ 
and $\xi^+(\tau) \in [b,\xi^+]$ such that
$$
v(\xi^-) +\tau= v(\xi^+)+\tau=v_{com}+\tau.
$$
Here, we change the notation and we write $\xi^-(\tau)$
(respectively, $\xi^+(\tau)$) in place of $a(\lambda)$ (respectively,
$b(\lambda)$ and $h = 1/\lambda$.

In order to solve (\ref{rn}), we have to find simultaneously $u$ and
$\si(x)\in \sgn u_x$, where $\sgn$ is understood as a maximal monotone
graph. We want that $u$ be constant equal to $v(\xi^-(h))$ on yet unspecified
$[\xi^-(h)$, $\xi^+(h)]$ containing $x_0$. On this interval $u_x$ will be zero
and $\si(x)\in \sgn 0$ will be different from zero.
Integration of  (\ref{rn}) over  $\xi^-(h)$, $\xi^+(h)$
yields an analogue of (\ref{na3}), i.e.
\begin{equation}\label{rnhta}
-2h = \int_{\xi^-(\tau(h))}^{\xi^+(\tau(h))}(u-v)\,dx.
\end{equation}
In Lemma \ref{podko} we established continuity of the mapping 
$$
[0,\tau_{max})\ni\tau\mapsto 
\int_{\xi^-(\tau)}^{\xi^+(\tau)}(v(\xi^-(\tau))-v(x))\,dx,
$$ 
(where $\tau_{max}=v_M-v_{com}$). Moreover, it is  strictly decreasing
and  equal to zero for $\tau=0$. 

Hence, for a fixed $h$ there is at most one $\tau(h)$ such that
(\ref{rnhta}). If there is such  $\tau(h)$, then for the sake of
simplicity we shall call $\xi^\pm(\tau(h))$ by $\xi^\pm(h)$. Thus, we
set
$$
u(x) =\left\{
\begin{array}{ll}
v(\xi^-(h)) & \hbox{for }x\in [ \xi^-(h),  \xi^+(h)],\\
v(x) &\hbox{elsewhere}.
\end{array}
\right.
$$
Of course,  we have the estimate
$\| u_x\|_p \le \|v_x\|_p $ for any $p\in[1,\infty]$.

We have to define $\si\in \sgn (u_x)$. 
On the set $\{ u_x>0 \}\cup \{
u_x<0 \}$, there is no problem for we put
$$
\si (x) = \sgn(u_x(x)).
$$

Before we proceed with the inductive  step we introduce a new
notation. Let us suppose that $F(a_1,b_1)$, $\ldots$, $F(a_N,b_N)$ are
all essential facets. Let us look at $[0,1]\setminus
\bigcup_{j=1}^N[a_i,b_i]$ consisting of open sets (in $[0,1]$)
$(p_j,q_j)$, $j=1,\ldots, N+2$. Each of the intervals $(p_j,q_j)$ has
the following property, either $v|_{(p_j,q_j)}$ is increasing, then
we write $(p_j,q_j)\in E^+(v)$, or $v|_{(p_j,q_j)}$ is decreasing, then
we write $(p_j,q_j)\in E^-(v)$. We note that the intervals $(p_j,q_j)$
are maximal sets (with respect to set inclusion) with the above
property. 

By the very definition, for $u$ as in the statement of this Lemma, we
have the following decomposition into disjoint sets,
\begin{equation}\label{dek}
[0,1] = E^+(u) \cup E^-(u) \cup \Xi_{ess}(u_x)\cup 
(\Xi(u_x)\setminus\Xi_{ess}(u_x).
\end{equation}
In general, if $u\in AF$ we say that $x_0\in E^+(u)$ (resp. $x_0\in
E^-(u)$), iff $x_0\not\in\Xi_{ess}(u_x)$ and there is
$(\alpha,\beta)$, a connected component of $\{ u_x>0 \}$, such that
there is $(l(\alpha),r(\beta))$ 
containing $(\alpha,\beta)$ and maximal with respect to set inclusion
such that $u|_{(l(\alpha),r(\beta))}$ is increasing. In a analogous
manner we define $E^-(u)$. We notice that $E^+(u)$ and $E^-(u)$ are
open and disjoint. We notice that $E^+(u)$ and $E^-(u)$ are
open and disjoint. It is obvious that the decomposition   (\ref{dek}) is valid for smooth $u$. Moreover, it is not difficult to notice (we will
not use it) that if  $u\in AF$, the decomposition   (\ref{dek}) holds.

We note that $\{ u_x>0 \}\subset E^+(v)$ and  $\{ u_x<0 \}\subset
E^-(v)$ with the possibility of strict inclusion. We set $\si$ equal
to 1 on $E^+(v)\setminus \{ u_x>0 \}$ and $\si$ equal
to $-1$ on $E^-(v)\setminus \{ u_x<0 \}$.

Otherwise we define $\si $ so that (\ref{rn}) holds, e.g. on
$[\xi^-(h),  \xi^+(h)]$ we set
$$
\si(x) =1 +\frac1h\int_{\xi^-(h)}^x(v(\xi^-(h)) -v(x))\,dx.
$$
The complement of 
$E^+(v)\cup E^-(v)\cup
[\xi^-(h), \xi^+(h)]$ is easy to consider and left to the reader.

We also mentioned the possibility that 
\begin{equation}\label{quat}
|\int_{\xi^-(\tau_{max})}^{\xi^+(\tau_{max})}(v(\xi^-(\tau))-v(x))\,dx|
=: 2h_{max}<2h.
\end{equation}
If this happens we proceed as follows.
We find $u$ by the above procedure yielding a minimizer of the
functional $J_{h_{max},v}$. By Lemma \ref{lac4}, we  split the
minimization problem into two: one for $J_{h_{max},v}$ already accomplished
and for $J_{h-h_{max},u}$. Let us notice that the process above for
$h=h_{max}$ yields $u$ which is monotone. We have already noticed
that if $u$ is monotone, then the unique minimizer of $J_{h-h_{max},u}$
is $u$ itself.  


Here  comes the inductive step. We construct $u$ for
$v$ with  $N+1$  essential facets, denoted as above,
provided
that we know how to deal with $v$ with $N$ essential facets.
For each essential facet $F(a_i,b_i)$, $i=1,\ldots, N+1$, we may find intervals
$[\xi^-_i,\xi^+_i]$,
$i=1,\ldots,N+1$, constructed as above. 
We may assume that the ordering is such that the
sequence of numbers
$\int_{\xi^-_i}^{\xi^+_i} |v(x)-v(\xi^-_i)|\,dx$, $i=1,\ldots,N$ is
decreasing. By the 
process described earlier, for a given positive $h,$ we define intervals
$[\xi^-_i(h),\xi^+_i(h)]$. We have two cases to consider: (a) interval
$[\xi^-_{N+1}(h),\xi^+_{N+1}(h)]$ is contained in $[0,1]$ and it does
  not overlap any of the intervals $[\xi^-_i(h),\xi^+_i(h)]$,
  $i=1,\ldots,N$, i.e.
$\xi^-_{N+1}(h)$ is positive,  and it is bigger than
$\xi^+_j(h)$ for all $j$ such that $\xi^-_{N+1} >\xi^-_j(h)$;
at the same time
$\xi^+_{N+1}(h)<1$ and $\xi^+_{N+1}<\xi^-_k(h)$
for all $k$ such that $\xi^+_{N+1}<\xi^+_k(h)$; (b) the
previous condition does not hold, i.e.  interval
$[\xi^-_{N+1}(h),\xi^+_{N+1}(h)]$ is not contained in $[0,1]$ or it
  intersect at least one interval $[\xi^-_i(h),\xi^+_i(h)]$.

The first case presents no problem. The intervals
$[0,\xi^-_{N+1}(h)]$, $[\xi^+_{N+1}(h),1]$ contain no more than $N$
essential facets $F(a_k,b_k)$. 
Thus, by the inductive assumption we know how to resolve any possible
overlapping. 

If (b) occurs, then there is $j_0$, such that
$[\xi^-_{j_0}(h),\xi^+_{j_0}(h)]$ intersects
$[\xi^-_{N+1}(h),\xi^+_{N+1}(h)]$ or $[\xi^-_{N+1}(h),\xi^+_{N+1}(h)]$
is not contained in $[0,1]$. 
The second case is easier, we shall deal with it first. It means that
there is $h_0<h$ such that $\xi^-_{N+1}(h_0)=0$ or
$\xi^+_{N+1}(h_0)=1$. But then, as we know, $F(0,\xi^+_{N+1}(h_0))$ or 
$F(\xi^-_{N+1}(h_0),1)$ are not essential facets, thus we consider the
minimization of $J_{h_0,v}$ with minimizer $u_0$ having $N$ essential
facets (of course we have to adjust the integral of integration in the
functional). If it is so, then by the inductive assumption we are able
to resolve 
any interactions, i.e. intersections of $N$ essential facets. Then,
we solve the minimization of $J_{h-h_0,u_0}$ where the minimizer has
no more than $N$ essential facets.

Thus, inevitably we  deal with interactions of
facets. Resolving the interactions is easier 
with Lemma \ref{lac4} below, which says that $h$ may be split, if
necessary, when $\xi^-_{j}(h_1)=\xi^+_i(h_1)$, and $h_1<h$,
while $\xi^-_{j}<\xi^+_i$. Let us assume that $h_1$ is the smallest
with this property. We solve our problem with $v$ and $h_1$, we 
find a minimizer of $J_{h_1,v}$. We may do
so, because of lack of interactions, we denote its solution by
$u^1$. Due to the occurrence of interactions the number of the
essential facets $F(a'_i,b'_i)$  of $u^1$ is smaller than for
$v$. Thus, we may use the inductive assumptions to continue, i.e. to
solve our problem with data $u^1$ and $h_2=h-h_1$, in place of $h$. By
Lemma \ref{lac4} solution $u^2$ is what we need. The proof of the
lemma is complete.
\qed

Our next
Lemma explains that $h$ may be split into smaller steps at will. This
permits to perform additional analysis at the intermediate steps.

\begin{lemma}\label{lac4}{\sl Let us suppose that $v$ is absolutely
    continuous and $h_1,$ $h_2>0$ the sets $\{v_x>0\}$, $\{v_x<0\}$ are
    open and they have a finite number of connected components. If
    $u^1$ is a minimizer of
$$
\cJ_{h_1,v}(u) =\int_0^1 h_1| u_x | +\frac 12 (u-v)^2
$$
while  $u^2$ is a minimizer of
$$
\cJ_{h_2,u^1}(u) =\int_0^1 h_2| u_x | +\frac 12 (u-u^1)^2,
$$
then $u^2$ is a minimizer of 
$$
\cJ_{h,v}(u) =\int_0^1 h| u_x | +\frac 12 (u-v)^2
$$
with $h=h_1+h_2$.}
\end{lemma}
{\it Proof.} In fact due to our
assumptions we have solutions to the equations
\begin{equation}
 h_1 \frac d{dx} \si^1 = u^1 - v, \qquad h_2 \frac d{dx} \si^2 = u^2 - u^1.
\end{equation}
We note that the sequence of implications: $u^2_x$ is different from zero at
$x$, then $u^1_x$ has a sign there, hence $v_x$  has a sign too. Moreover, if
$v_x=0$ on an interval $(\al,\be)$, then $u^1_x$,  $u^2_x$ are zero
$(\al,\be)$ too. 

We want to show that 
\begin{equation}\label{rn10}
h\frac d{dx}\sgn u^2_x = u^2 -v
\end{equation}
has a solution. Let us add up the two equations above. This yields,
$$
h\frac d{dx}\left(\frac {h_1}h\si^1 + \frac {h_2}h\si^2  \right)
= u^2- v.
$$
Of course $\si := \frac {h_1}h\si^1 + \frac {h_2}h\si^2 \in[-1,1]$. If at $x$
we have $u^2_x(x)>0$, then $v_x(x)>0$. Hence,
$$
 \si(x) = \frac {h_1}h\si^1(x) + \frac {h_2}h\si^2(x) = 
\frac {h_1}h +\frac {h_2}h =1.
$$
The situation is similar if  $u^2_x(x)<0$.
Let us suppose now that  $u^2_x(x)=0$, then regardless of the sign of 
$u^1_x(x)$, we know that $\si (x) \in[-1,1]$ and, by the definition
of $\si$, equation (\ref{rn10}) is satisfied. In particular,
$$
-2h = h\int_{\xi_2^-}^{\xi_2^+}\si(x)\,dx
=\int_{\xi_2^-}^{\xi_2^+}(v(\xi_2^-)-v(x))\,dx=
\int_{\xi_2^-}^{\xi_2^+}(u^2(\xi_2^-) -v(x))\,dx . 
$$
\qed

The value of this result is that it permits us to split $h$. We may say
that this shows the semigroup property. 
Finally, we show that functions with finite number of essential
facets are dense in the topology of $W^1_2$.

\begin{lemma}\label{lac2}{\sl If $v$ is smooth with $v(a)=A$,
    $v(b)=B$, then there exist $v_k$ satisfying the assumption of Lemma
    \ref{lac1}. Moreover $v_k$ converges weakly to $v$ in  $W^{1}_{2}$ and
    $\| v_k\|_{1,2}\le \|v\|_{1,2}$.}
\end{lemma}
{\it Proof. } The sets $E^+(v)$, $E^-(v)$ consist of at most
countable number of open intervals,
$$
E^\pm(v) = \bigcup_{k\in{\cal I}} I_k^\pm(v). 
$$
Subsequently, we suppress the $\pm$ superscripts.

We order the intervals $I_k$, $k\in\bN$ in  so that $|I_k|\ge
|I_{k+1}|$. On $ \bigcup_{j=1}^k I_j$,
we set $v^k(x)=v(x)$. On the complement, we define $v^k$ to be
piecewise linear and continuous. We immediately notice that
$$
\|v^k_x\|_2\le \|v_x\|_2,
$$
because the linear functions are harmonic. Hence, they minimize the
functional $\int |v_x|^2$ with Dirichlet data. We have to show that
$v^k_x$ converges to $v_x$ in $L_2$. 

We will show first the pointwise convergence of $v^k$. Let us take any
$x\in[0,1]$. If $x\in E^+(v)\cup E^-(v) $, then $x\in I_{j_0}$, hence $v^k(x)=
v(x)$ for $k\ge j_0$. We suppose now that $x_0$ is in the complement
of $E^+(v)\cup E^-(v) $. For the sake
of further analysis, we set $\cF_k =[0,1]\setminus\bigcup_{i=0}^k
I_i$. Each of the sets $\cF_k$ consists of a finite sum of closed
intervals and $x_0\in[\al_k,\be_k]$, $k\in\bN$.
By  construction the sequence $\al_k$ is
increasing, while $\be_k$ is decreasing. We shall call by $\al$ and
$\be$ their respective limits. Of course, we have that
$v^k(\al_k)=v(\al_k)$ thus this sequence convergence to $v(\al)$, while
$v^k(\be_k)=v(\be_k)$ converges to $v(\be)$. We have two case to
consider: 1) $\al < \be$, 2) $\al = \be$. In the first case we have
$\displaystyle v^k_x =\frac{v(\be_k)-v(\al_k)}{\be_k-\al_k}$. This
must converge to 
zero. Otherwise, we had $v_x\neq 0$ on a subset of $(\al,\be)$ of
positive measure which is impossible by the  definition of
$\cF_k$'s. Hence, $v(\al)=v(x)=v(\be)$ for all $x\in[\al,\be]$,
i.e. $v^k(x)$ converges to $v(x)$. 

If $\al = \be$, then our reasoning is similar and by continuity of
$v$ we deduce that $v(\al)=v(x)=v(\be)$.

Thus, we have shown that $v^k$ converges everywhere to $v$. On the
other, hand the bound $\|v^k_x\|_2\le \|v_x\|_2$ implies that we can
select a weakly convergent subsequence. Due to uniqueness of the limit
it must be $v$. Since any convergent subsequence of $v^k$ converges
to $u$, the whole sequence  $v^k$ converges to $v$.

Moreover,  due to Sobolev embedding, we deduce that $v^k$ converges to
$v$ uniformly.
\qed

\begin{lemma}\label{lac3}{\sl If a sequence  $v^k\in W^{1}_{2}$ converges
    to $v$ in  $W^{1}_{2}$ and  $u_k\in W^{1}_{2}$ is the  sequence of
    corresponding minimizers of $\cJ_{v^k}$, then $u^k$ converges to $u$
    weakly in  $W^{1}_{2}$ and strongly in $L_2$. Moreover,  $u$ is a
    minimizer of $\cJ_{v}$ and $\| u\|_{1, 2}\le \|v\|_{1, 2}$.
}
\end{lemma}
{\it Proof. } The convergence of  $u^k$ in  $L_2$ follows from the
monotonicity of the subdifferential. Indeed, since $u^k$ is a
minimizer, then
$$
h\partial\cJ(u^k)+ u^k -v^k\ni 0,
$$
i.e. there exists $\zeta_k\in\partial\cJ(u^k)$ such that for any test
function $\phi\in L_2$ we have,
$$
h\langle \zeta_k,\phi\rangle + \langle u^k,\phi\rangle 
= \langle v^k,\phi\rangle .
$$
Once we take $\phi = u_k- u_l$, we can see that
$$
h\langle \zeta_k- \zeta_l, u^k- u^l\rangle  + \|u^k- u^l\|_2^2
= \langle v^k- v^l, u^k- u^l\rangle .
$$
Due to monotonicity of $\partial\cJ$ this implies that $\|u^k- u^l\|_2
\le \|v^k- v^l\|_2 $. Thus, $u^k$ 
converges in $L_2$ to $u^*$.

The estimates we have already shown yield
$$
\|u^k_x\|_p\le \|v^k_x\|_p\le \|v_x\|_p +1,
$$
for sufficiently large $k$. It means, that we can select a weakly convergent
subsequence with limit $\bar u$. Due to uniqueness of the limit,
$u^*=\bar u$. Moreover, all weakly convergent subsequences have a common
limit $u^*$. Hence the sequence $u^k$ converges weakly in $W^{1}_{2}$ to
$u^*$. 

We know that  $\cJ_v$ has a unique minimizer $u$. Now, we have to show
that $u^*$ is the minimizer of $\cJ_v$, i.e. $u^*= u$. Obviously, we have
\begin{equation}\label{rnkn}
J_{v^k}(u) \ge J_{v^k}(u^k).
\end{equation}
Due to the lower semicontinuity of the $BV$ norm, we have
$$
\liminf_{k\to\infty} J_{v^k}(u^k) =\liminf_{k\to\infty}
h\int_a^b |u^k_{x}| +\lim_{k\to\infty}\frac12\int_a^b(u^k-v^k)^2 
\ge h\int_a^b |u^*_{x}| +\frac12\int_a^b(u^*-v)^2 
=J_{v}(u^*).
$$
On the other hand, we have
$$
\lim_{k\to\infty} J_{v^k}(u) =
\int_a^b h |u_{x}| +\lim_{k\to\infty} \frac12\int_a^b (u-v^k)^2 
=J_v(u).
$$
Thus, due to (\ref{rnkn}), we infer that
$$
J_v(u) \ge J_{v}(u^*).
$$
Since $u$ is a unique minimizer of $J_{v}$, we conclude that $u=u^*$.
 Our claim follows.
\qed

\bigskip
We are ready to show our main results.

{\it Proof of Theorem \ref{twac}}. 
{\it Step 1.} We have already  noticed in Lemma \ref{bla} that there
exists a minimizer of $J_{h,v}$. Hence, there exists a solution to
the following inclusion
$$
h\partial\cJ(u)+ u-v \ni 0.
$$
Moreover, it is also unique, because if we had two, say $u^1$ and $u^2$, then
for some $\zeta^i\in\partial\cJ(u^i)$, $i=1,2$, we had
$$
u^2-u^1+ h(\zeta^2-\zeta^1 )=0.
$$
Once we apply the test function $u^2-u^1$ to both sides, we see that
$\|u^2-u^1\|_2 \le 0.$
Hence the  claim, i.e. for any $v\in L_2$ there exists $u\in D(\cJ)$ a unique
minimizer of $\cJ_v$. 
The above argument yields only that $u$ belongs to $ BV$. 
Now, the goal is to improve regularity of minimizers.

{\it Step 2.} We will call by $\bar v_\ep$ the standard mollification of $v$.
Of course, $ \|\bar v_\ep\|_{1,p} \le \|v\|_{1,p}$, but $\bar v_\ep$
may not satisfy the boundary conditions, so we add a linear function. We
call the result by $ v_\ep$. Of course, 
$ \| v_\ep\|_{1,p} \le  \|\bar v_\ep\|_{1,p} + O(\ep)$.

We will show that the sequence of solutions $u_\ep$ to the minimization 
problem converges weakly in $W^{1}_{p}$ and strongly in $L_2$ to $u$ a solution
to the original problem.

{\it Step 3.} Since $v_\ep$ is smooth, then the sets
$E^+(v_\ep)$, 
$E^-(v_\ep)$ 
which we defined in Lemma \ref{lac1} 
are open, i.e.
$$
E^+(v_\ep)=\bigcup_{i\in I^+}(\al^+_i,\be^+_i),\quad
E^-(v_\ep)=\bigcup_{i\in I^-}(\al^-_i,\be^-_i).
$$
The index sets $I^+$, $I^-$ are at most countable. We may arrange the intervals
at will.

{\it Step 4.} We know by Lemma \ref{lac1} above, that if $v$ is smooth
and the sets 
$I^+$, $I^-$ are finite, then $u\in W^{1}_{p}$, for any $p\in[1,+\infty)$
and it is piece-wise smooth. 
Moreover,
$v=u$ on $E^+(u)\cup E^-(u)$. In particular, the set 
$[a,b]\setminus E^+(u)\cup E^-(u)$ is a finite sum of closed
intervals, so that we may write 
$$
[a,b]\setminus E^+(u)\cup E^-(u)=\bigcup_{i=1}^N[\xi^-_i,\xi_i^+].
$$
In particular, it is possible that $\xi^-_i=\xi_i^+$.

We also know that if for some $\de>0$ function $v$ is monotone on
$[\xi^-_i-\de,\xi_i^+ +\de]$, then $u=v$ on
$[\xi^-_i,\xi_i^+]$, i.e. $v([\xi^-_i,\xi_i^+])$ is a zero curvature
facet. More interesting is the case, when for some $\de>0$ function
$v$ is convex or concave on $[\xi^-_i-\de,\xi_i^+ +\de]$. Then, 
$u = v(\xi^-_i) = v(\xi^+_i)$ on $[\xi^-_i,\xi_i^+]$ and
$$
\int_{\xi^-_i}^{\xi^+_i}(u(x)-v(x))\,dx = 2h.
$$
From these properties, we deduce that 
\begin{equation}\label{szawp}
\| u_x\|_p\le\|v_x\|_p,
\end{equation}
for all $p\in[1,\infty)$.

{\it Step 5.} In Lemma \ref{lac2}, we constructed a sequence of
continuous, piecewise 
smooth $v_\ep^k$ converging weakly to $v_\ep$ in $W^1_2$.

Let us call by   $u^k_\ep$ the minimizers of $\cJ_{v^k_\ep}$.
Monotonicity of 
$\partial\cJ$ implies convergence of $u^k_\ep$ in $L_2$. Indeed,
if $u^k_\ep + \partial\cJ(u^k_\ep) \ni v^k_\ep$, then taking difference 
and applying it to a test vector yields,
$$
(u^k_\ep - u^l_\ep , p) +h\langle j^k -j^l, p\rangle =
(v^k_\ep - v^l_\ep , p),
$$
where $j^k\in\partial\cJ(u^k_\ep)$, $j^l\in\partial\cJ(u^l_\ep)$.
When we choose $p= u^k_\ep - u^l_\ep$, then monotonicity of the 
subdifferential implies
$$
\|u^k_\ep - u^l_\ep\|_2 \le \|v^k_\ep - v^l_\ep\|_2.
$$
Hence, the $L_2$  convergence of $v^k_\ep$ implies  the $L_2$  
convergence of $u^k_\ep$ to a limit $u_\ep$. We have to improve the 
regularity of the limit. For this purpose, we notice that the
estimate (\ref{szawp}) applied to the 
sequence $v_\ep^k$ yields,
$$
\| u^k_{\ep,x}\|_p\le\|v^k_{\ep,x}\|_p
$$
for any $p\in(1,2]$. Hence, we can select a weakly convergent 
subsequence in $W^{1}_{2}$ with limit $u_\ep^\infty$. Due to uniqueness 
of the limit we conclude that $u_\ep=u_\ep^\infty$, i.e. $u_\ep$ is
in $W^{1}_{2}$ for any finite $p$. This also implies that $u^k_\ep$ 
converges to $u_\ep$ uniformly.

Since the norm is lower semicontinuous 
we also infer that
$$
\| u_{\ep,x}\|_p\le\|v_{\ep,x}\|_p\le \|v_{x}\|_p.
$$
So the same argument permits us to pass to the limit with $\ep\to 0$
to conclude that $u_\ep$ converges to a limit $u$ strongly in $L_2$,
$L_\infty$ and weakly in $W^1_2$.

{\it Step 6.} We have to show that $u_\ep$, for $\ep>0$, and $u$ are
minimizers of $\cJ_{v_\ep}$ for the corresponding data $v_\ep$ or
$v$. For this purpose, we invoke  Lemma \ref{lac3}.
\qed

\bigskip
We also note a conclusion from the proof of Theorem \ref{twac}.
\begin{corollary}{\sl Let us suppose that 
$v$ is continuous and piecewise smooth, such that
    one sided derivatives exit everywhere. The sets $\{v_x>0\}$,
    $\{v_x<0\}$ are open with a finite number of connected components
    denoted by $K$.
Then, $u$ the unique minimizer of $\cJ_v$,  belongs to
$W^1_p$, for any $p\in[1,+\infty)$ and it is piecewise smooth. 
Moreover,
$v=u$ on $E^+(u)\cup E^-(u)$ and there exists $\si\in W^1_1$, such that
$\si(x)\in \sgn u_x(x)$ and
$$
-h\frac d{dx}\si = v-u.
$$
Furthermore,
$\| u_x\|_p\le\|v_x\|_p,$
for all $p\in[1,\infty)$.}

\end{corollary}

Theorem \ref{twac} is slightly too general for our purposes, Theorem
\ref{wngl} is its  refinement. We will prove it momentarily.

\smallskip

{\it Proof of  Theorem \ref{wngl}.}
Part (a) is obvious when $K(v)=\infty$. If $K(v)<\infty$, then the claim
follows from the construction of $u$ if $h$ is sufficiently small. For
a general $h$ we have to use Lemma \ref{lac4}. 

Our proof of part (b) starts with the observation that $v_x\in BV$
implies  $v_x\in L_\infty$. Hence, we can pass to
the limit with $p$ in the estimate
$\|u\|_{1,p}\le \|v\|_{1,p}.$
Thus, $\|u\|_{1} +\|u_x\|_{1} \le \|v\|_{1}+ \|v_x\|_{1}$. 

If $v_x\in BV$, then by the general theory, see e.g. \cite{ziemer},
there exists a sequence of smooth functions, $\{v_k\}$, such that
$\|v_{k,x}\|_{BV}$ converges to $\|v_{x}\|_{BV}$. We apply Lemma
\ref{lac3} to deduce existence of a sequence  $\{v_{km}\}$ such that
the sets $\{v_{km,x}>0\}$ and  $\{v_{km,x}<0\}$  are open and have a
finite number of components. Moreover,
$\displaystyle{\lim_{m\to\infty}v_{km}=v_k}$ in $W^{1}_{1}$.

Now, it is easy to calculate  the norm $\|u_{km,x}\|_{BV}$ for the
corresponding minimizers $u_{km}$ for sufficiently small $h$. We have 
\begin{eqnarray*}
\int_a^b |D u_{km,x}| &=& 
\sum_{i} \int_{(\xi_i^+(h),\xi_{i+1}^-(h))} |D v_{km,x}| 
+ \sum_{i}( |v_{km,x}^+(\xi_i^+(h))|+|v_{km,x}^-(\xi_{i+1}^-(h))|)\\
 &\le&\sum_{i}\int_{(\xi_i^+(h),\xi_{i+1}^-(h))} |D v_{km,x}| \\&&+
\sum_{i}( |v_{km,x}^+(\xi_i^+(h))-v_{km,x}^+(\xi_i^+(0))|
+|v_{km,x}^-(\xi_{i+1}^-(h))-v_{km,x}^-(\xi_{i+1}^-(0))|) \\&&+
\sum_{i}( |v_{km,x}^+(\xi_i^+(0))|+|v_{km,x}^-(\xi_{i+1}^-(0))|)
\\
&\le&  \sum_{i}\int_{(\xi_i^+(0),\xi_{i+1}^-(0))} |Dv_{km,x}|+
\sum_{i}( |v_{km,x}^+(\xi_i^+(0))|+|v_{km,x}^-(\xi_{i+1}^-(0))|) \\
 &=& \int_a^b |D v_{km,x}|.
\end{eqnarray*} 
Here, we use the convention that if $\xi^+_1(h)>a$, then we write
$\xi^+_0(h)=a$ and $\xi^-_{N+1}=b$ provided that  $\xi^-_{N}<b$.

That is, we have 
$$
\|u_{km,x}\|_{BV} \le \|v_{km,x}\|_{BV} .
$$
We can find  $m_k$ converging  to zero as $k$ goes to infinity such that
$\|v_{km_k,x}\|_{BV} \le \|v_{k,x}\|_{BV}+1/k.$
Finally, we use \cite[Theorem 5.2.1]{ziemer} to conclude  that
$$
\|D u_x\| \le \liminf_{k\to\infty}\|D u_{km_k,x}\|\le
\lim_{k\to\infty} (\|D v_{k,x}\|+1/k)=\|v_{x}\| .\eqno\Box
$$


\section{Asymptotics and  examples} 

Here, we present the proof of Theorem  \ref{asymp}, an example of an
explicit solution  and numerical results describing the time behavior
of solutions.

\subsection{A proof of Theorem \ref{asymp}.}
Here is the argument.
There is a finite number $N$
of facet merging events 
$$
0=t_0<t_1< \ldots < t_N< \infty,
$$
when $u$ has no time derivative but only the right-time derivative. Moreover, $N\le K_{ess}(u_{0,x})$. We shall estimate $\max_{i=0,\ldots N-1} \{ t_{i+1} -t_i\}$. Let us set
$$
B= \max\{ a_b, a_e\},\quad b= \min\{ a_b, a_e\},\quad
\Delta_M = \max u_0(x) - B,\quad \Delta_m = b - \min u_0(x),
$$ 
and $\ell = 1$ is the length of $I=[0,1]$. We notice that since our solution is almost classical, then $u_t$ exists except $t\in\{t_0,t_1, \ldots,  t_N\}$.
Moreover, $u_t$ is the vertical velocity of $u$. It is obvious from the definition of the composition $\bar\circ$ that the absolute value of $(\hbox{sign}\,\bar\circ u_x)_x$ is bigger or equal $2/\ell$. We notice that the distance each essential facet travels in the vertical motion between collisions is no bigger than
$$
A = \max \{ \Delta_M, \Delta_m, B-b\}.
$$
Since we have a lower bound on the vertical velocity of $u$, we conclude that 
$$
\max_{i=0,\ldots N-1} \{ t_{i+1} -t_i\}\le A\cdot \frac 2\ell.
$$
Thus, we have the following estimate 
\begin{equation}\label{exte}
t_{ext} \le 2 K_{ess}(u_{0,x}) A/ \ell.
\end{equation}
Hence  $K_{ess}(u_x(t_{ext})=0$, then thus $u(t)$ for $t\geq t_{ext}$
is a monotone function being a stationary state of the system. 
\qed

\subsection{An explicit solution}
In order to illustrate the behavior of a particular solution we take
$x^2$ as an  initial datum for (\ref{i1}). We consider this system on the
interval $(-1,1)$,
\begin{equation}\label{example1}
 \begin{array}{lcr}
 u_t - \frac{d}{dx} \sgn u_x =0 & \mbox{ in } & (-1,1)\times (0,T), \\
u(-1,t)=u(1,t)=1 & \mbox{ for } & t \in (0,T),\\
u|_{t=0} = x^2 & \mbox{ for } & (-1,1).
\end{array}
\end{equation}
The proved results quarantee us the following form of the solution to 
(\ref{example1}),
\begin{equation}\label{example2}
 u(x,t)=\left\{
\begin{array}{lcr}
a^2(t) & \mbox{ for } & |x| \leq a(t), \\
x^2 & \mbox{ for } & |x| \in (a(t),1)
\end{array} \right.
\end{equation}
By Definition \ref{ddxx} we get that 
$$
\frac{d}{dx} \sgn \bar\circ u_x |_{[-a(t),a(t)]} = \frac{1}{a(t)}.
$$ 
Thus by (\ref{example1}) and (\ref{example2}) we find a relation on
$a(t)$ 
as follows 
$$
\d_t a^2(t)=\frac{1}{a(t)}, \mbox{~~~ hence ~ } a(t)=\sqrt[3]{\frac 32 t}
$$
to keep the agreement to the initial datum.

Summing up the length of the facet is $2a(t)=2\sqrt[3]{\frac 32 t}$,
the speed of it is $\d_t a(t) \sim t^{-2/3}$ and
the extinction time of 
$u\equiv 1$ is $T_{stab} = \frac 23.$

\subsection{Numerical simulations}
Now, we are prepared to computer implementations of our results.
Simulations were done in Octave package. The main part of the program
is a loop running until the graph reaches it's final shape. During one
step all facets (i.e. points where $0 \in \partial f$) are moved until
(if it is possible) each of them fills the area equal to $2h$. In the
pictures shown below we used $h=5$. The reason why it may be not
possible to fill the $2h$ area is that the moving facet may reach the
boundary of the interval that it is defined on or it may reach the
boundary of another facet after it filled the required area (whereas
each of them moved separately may fit its domain). When any of these
interactions happens, we change the $h$ value for a maximum reached
value (let us call this new value $h_{min}$) and move all facets so
that they fill the area of $2h_{min}$. We use $h_{min}$ just in this
one step but for all facets and then get back to $h$ value. After each
step, we recalculate domains and check if we still use all functions
(some of them may disappear, as the $x^2-2x$ function defined on
$[0,1]$ interval after the first step of the $v_1$ example from table
1). 

In none of the presented examples 
a facet fills the maximum area. We chose $h$ big enough to avoid unnecessary steps. 

We calculate the time a step takes as $\frac{2h_{min}}{2h}$. We do this using the following logic --- we make an assumption that one full step (i.e. area of $2h$ is filled) is my time unit, two full steps count as $t=2$, $\frac{1}{3}h$ takes $t=\frac{1}{3}$ to fill. In the pictures accumulated time is presented.

As an initial data in three presented examples, we use functions
described in the table below. The first column contains intervals
which set the domain, the next three columns contain formulas for
respective examples:
\begin{center}
\begin{tabular} {|l||l|l|l|}
\hline domain    & $v_1$  & $v_2$ & $v_3$ \\ \hline
$\left[-1.5,-1\right] $ & $ x-2       $ & $ 3x^2+11x  $ & $ 3x^2+11x $ \\ \hline
$\left[-1,0\right]    $ & $ -x^2+x+2  $ & $ -x^2+5x+1 $ & $ -x^2+5x+1 $ \\ \hline 
$\left[0,1\right]     $ & $ x^2-2x    $ & $ x-2       $ & $ 0 $ \\ \hline 
$\left[1,2\right]     $ & $ -x^2+5x+1 $ & $ 2x-7      $ & $ 2x-7 $ \\ \hline 
$\left[2,3\right]     $ & $ x^2-6x+8  $ & $ x^2-6x+8  $ & $ 1 $ \\ \hline 
$\left[3,4\right]     $ & $ 0         $ & $ -x^2+x+2  $ & $ -x^2+x+2 $ \\ \hline 
$\left[4,5\right]     $ & $ 2x-7      $ & $ x^2-2x    $ & $ x-2 $ \\ \hline 
$\left[5,5.5\right]   $ & $ 1         $ & $ x^2+15x   $ & $ x^2+15x $ \\ \hline
\end{tabular} \\
Table 1. Examples 1, 2, 3 (respectively) used in the simulations
\end{center}

To create the three examples, we use the same domain and permute functions to obtain interesting shape. In some cases, we have to move parts defined on some intervals vertically to obtain continuous result. Therefore, in some cases the same function used on the same interval has different values. What is more, we move the whole graph vertically so that the smallest value is 1; it makes integration easier without changing the shape of solutions.

We use polynomials as an approximation of a continuous function defined
on closed interval; in the examples mentioned they are of degree 2,
but the algorithm remains the same for polynomials of higher
degree. Functions defined on intervals model situation of
non-continuous derivative.

Let us look at  results of simulations presented on the figures:

\begin{figure}[!htb]
\centering
{\epsfxsize=12cm \epsfbox{
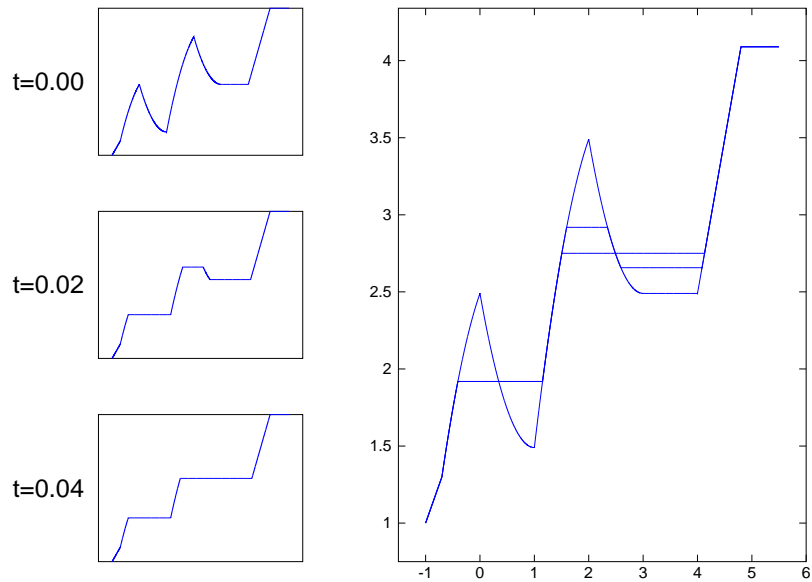}\\
\caption {The first example}}
\end{figure}

\begin{figure}[!htb]
\centering
\epsfxsize=12cm \epsfbox{
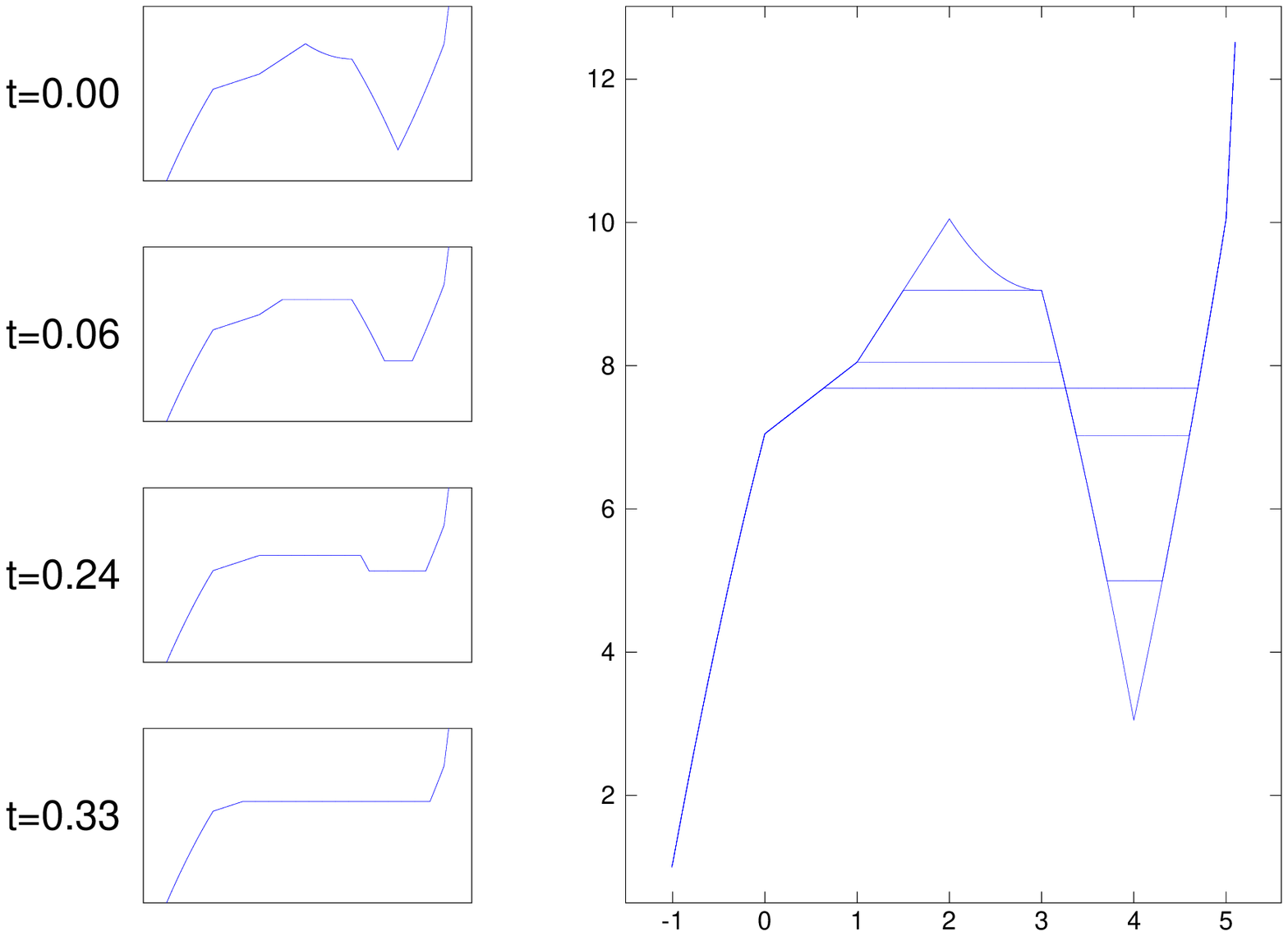}\\
\caption {The second example}
\end{figure}

\begin{figure}[!htb]
\centering
\epsfxsize=12cm \epsfbox{
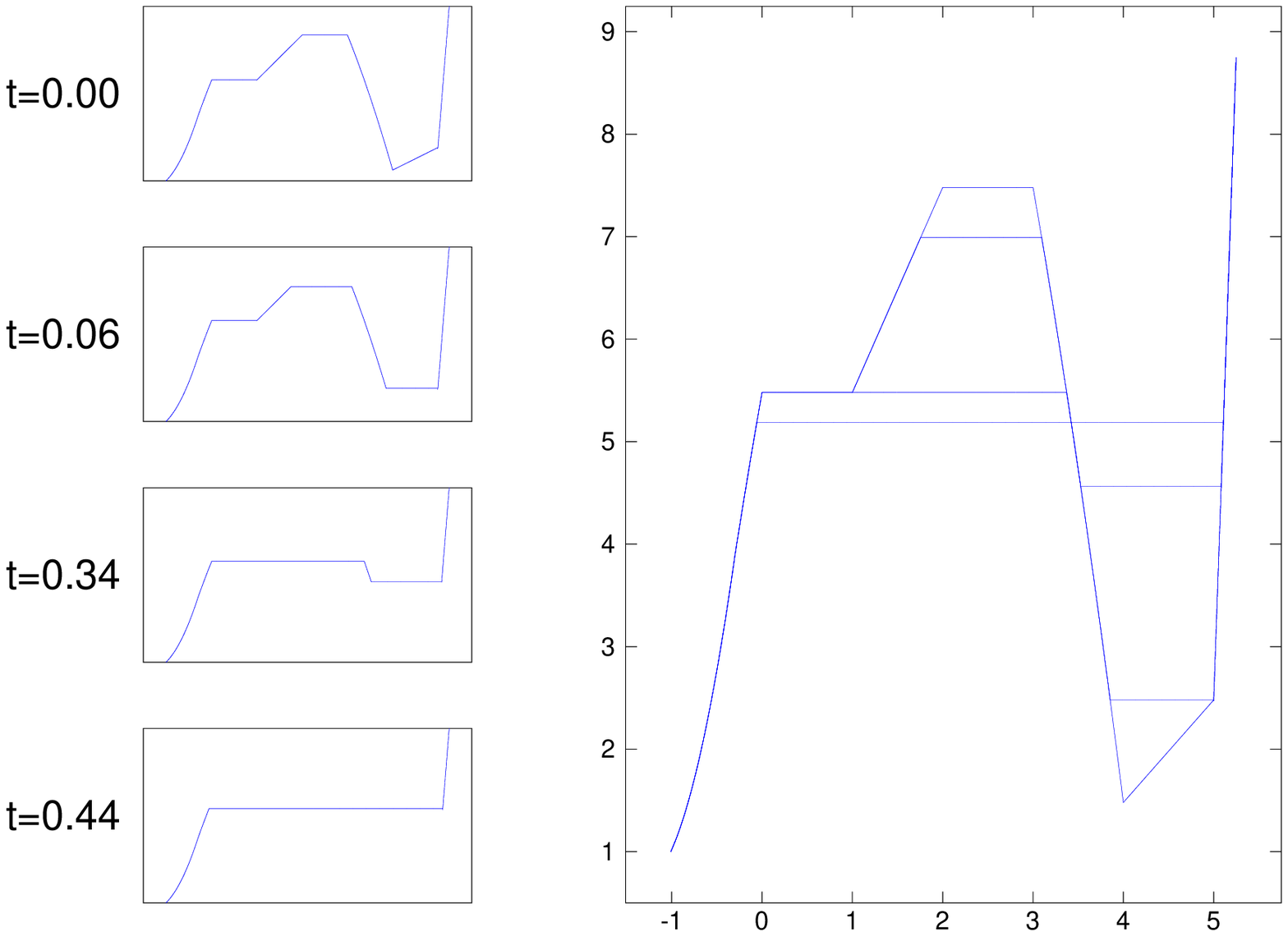}\\
\caption {The third example}
\end{figure}

Observe that, all degenerated facets disappear after
the first step of evolution. The number of regular facets that may
appear is limited by their number and the overall number of regular
facets decreases from the second step of evolution. The flat area
broadens with each step. All solutions remain continuous and their
$||\cdot||_{L_\infty}$ norm is bounded by the norm of initial data.

\bigskip\noindent{\bf Acknowledgment} PR thanks Professor Jos\'e
Maz\'on for inspiring conversations on the topic of this paper which
lead to an improvement of the  proof of Theorem
\ref{twac}. Special thanks go to the Iberia airline for creating extra
opportunities to work on this paper and on related topics.
The  work has been partly supported by
MN grant No. N N201 268935.

\end{document}